\documentclass{elsart}
\usepackage{etex}

\usepackage[numbers]{natbib}
\usepackage{easymat}
\usepackage{easybmat}
\usepackage{mathtools}
\usepackage{pslatex,times}
\usepackage{amssymb}
\usepackage{mathpazo}%psfig,epsfig,subfigure
\usepackage{easyeqn}
\usepackage[leqno,fleqn,intlimits]{empheq}

\usepackage[T1]{fontenc}
\usepackage[isolatin]{inputenc}
\usepackage{graphicx}
\usepackage{psfrag}
\usepackage{wrapfig}
\usepackage{multirow,tabls}
\usepackage{booktabs}

\usepackage{caption}

\usepackage[a4paper,textwidth=16cm, textheight=23cm,
					  hmarginratio=1:1,vmarginratio=1:1]{geometry}
\newcommand{\dt}{\delta t}
\newcommand{\pd}{\partial}

\newcommand{\bu}{{\mathbf u}}
\newcommand{\bom}{\boldsymbol{\omega}}

\newcommand{\lap}{\Delta}
\newcommand{\source}{S}
\newcommand{\bsource}{{\mathbf S}}

\newcommand{\bnabla}{\boldsymbol{\nabla}} 

\newcommand{\grad}{{\nabla}}

\newcommand{\curl}{{\nabla} \times}

\newcommand{\er}{\hat{\mathbf e}_r}

\newcommand{\ez}{\hat{\mathbf e}_z}
\newcommand{\et}{\hat{\mathbf e}_\theta}

\newcommand{\halfM}{\left\lfloor\frac{M}{2}\right\rfloor}

\newcommand{\pdt}{\partial_t}
\newcommand{\pdr}{\partial_r}
\newcommand{\pdrz}{\partial_{rz}}
\newcommand{\pdth}{\partial_\theta}
\newcommand{\pdthz}{\partial_{\theta z}}
\newcommand{\pdz}{\partial_z}

\newcommand{\pmh}{\pm\frac{h}{2}}

\newcommand{\mT}{\mathcal{T}}

\newcommand{\mQ}{\mathcal{Q}}
\newcommand{\ir}{\frac{1}{r}}

\newcommand{\iR}{Re^{-1}}

\newcommand{\laph}{\lap_h}

\newcommand{\ie}{i.e.}

\newcommand{\angvel}{\omega}

\numberwithin{equation}{section}
\begin{document}
\begin{frontmatter}
\title{
Poloidal-toroidal decomposition in a finite cylinder.\\
II. Discretization, regularization and validation}
\author{Piotr Boronski and Laurette S. Tuckerman}
\address{LIMSI-CNRS, BP 133, 91403 Orsay, France}

\begin{abstract}
The Navier-Stokes equations in a finite cylinder are written in terms of poloidal and toroidal
potentials in order to impose incompressibility.
Regularity of the solutions is ensured in several ways:
First, the potentials are represented using a spectral basis
which is analytic at the cylindrical axis.
Second, the non-physical
discontinuous boundary conditions at the cylindrical corners are smoothed
using a polynomial approximation to a steep exponential profile.
Third, the nonlinear term is evaluated in such a way as to eliminate singularities.
The resulting pseudo-spectral code is tested using exact polynomial solutions and the spectral convergence
of the coefficients is demonstrated.
Our solutions are shown to agree with exact polynomial solutions and with previous axisymmetric calculations of vortex breakdown and of nonaxisymmetric calculations of onset of helical spirals.
Parallelization by azimuthal wavenumber is shown to be highly effective.
\end{abstract}
\begin{keyword}
cylindrical coordinates, coordinate singularities, pseudo-spectral, recursion relations, radial basis, vortex breakdown, spectral convergence, poloidal-toroidal
\end{keyword}
\end{frontmatter}

\section{Introduction}
\label{sec:intro}

The von K\'arm\'an flow owes its name to T. von K\'arm\'an \cite{vonKarman}
who in 1921 first studied the flow in the semi-infinite domain bounded by a
single rotating disk using a similarity transformation.  In 1951,
\citet{Batchelor} extended the problem to the flow confined between two
infinite rotating disks. For rotating disks of finite radius, the
configuration can be described by three control parameters: the ratio $s$ of
the angular velocity of the two disks, the height-to-radius ratio $h$ and a
Reynolds number $Re$ based on the radius and the azimuthal velocity of one of
the disks.  The variation of these three parameters $(Re,s,h)$ has proved to
yield a rich variety of qualitatively different accessible flows, even before
the onset of turbulence.  The symmetries influence the transitions that the
flow can undergo.    This configuration is extensively studied in the
context of transition to complex and turbulent flows. All of these properties
explain why the von K\'arm\'an flow is increasingly considered as one of the
classical hydrodynamic configurations and why the scientific community is
interested in further exploring its complex behavior.

The first numerical studies were necessarily devoted to axisymmetric flows and
their stability.  In the rotor-stator configuration ($s=0$), vortex breakdown
forming characteristic recirculation bubbles was observed by number of authors
(\citet{Lugt}, \citet{DaubeSorensen89}, \citet{Lopez90}, \citet{Daube92}).
This now well-documented configuration became a benchmark for testing
axisymmetric codes.  Following \citet{LopezShen1998}, \citet{Speetjens} and a
number of other authors, we will validate our method in the axisymmetric
configuration by reproducing the stationary state at $Re=1800$ and, for
$Re=2800$, the non-stationary, oscillating flow, for which we will compare the
bifurcation threshold and the oscillation frequency against previous findings.
%(CITE OTHER AXI LARGE-ASPECT-RATIO ROTOR-STATOR WORK BY DAUBE/LE QUERE/COUSIN
%FROM LIMSI?)

Increasing computational power has made it possible to
study three-dimensional instabilities.  Breaking of axisymmetry has been the
subject of several studies, of which we mention these of \citet{Gauthier99},
\citet{Gelfgat}, \citet{Blackburn}, \citet{Lopez2001} and \citet{Nore04}.
Three-dimensional instability precedes axisymmetric instability for $h<1.6$
and $h>2.8$.  As the test problem for validating our code in three dimensions,
we have selected a configuration with $(s=-1,h=3.5,Re=2150)$, where the
perturbation of the axisymmetric flow takes the form of a helical spiral. In
section \ref{sec:test3d} we compare our results with those of
\citet{LopezMarquesShen2001} and \citet{Gelfgat}.

Interesting phenomena can also be observed in the turbulent regime.  
%According
%to the Taylor-Proudman theorem, sufficiently rapid rotation causes a flow to
%be independent of the direction of the rotation axis. In this regime one can
%expect to observe two-dimensional turbulence. This was explored by Swinney et
%al. \cite{Baroud,Jung}.  A related problem concerns laws describing decaying
%two-dimensional turbulence, which has been investigated mainly in Cartesian
%geometry (e.g. \citet{Yin}). (REALLY?) However, decaying two-dimensional
%turbulence in a system with $SO(2)$ symmetry is still not very well
%understood. Some results have been provided by \citet{Leprovost}.  (DOES THIS
%FIT?) (REF TO MOISY).  
Turbulence and large-scale structures may coexist.
In experiments in a turbulent counter-rotating configuration, \citet{Marie}
and \citet{Ravelet} discovered that a two-cell mean flow with a shear layer at
the cylinder mid-plane undergoes switching to a one-cell mean flow whose
shear layer is adjacent to the less rapidly rotating disk.  This transition
can be observed at Reynolds numbers which are numerically accessible.
%$Re\approx 5000$.
						
A comprehensive classification of the solutions for different values of
the parameters $(s,h,Re)$ is beyond the scope of this work.  Our main
purpose here is to develop a mathematical and algorithmic tool which can be
applied to von K\'arm\'an flow and to rotating turbulence, and which can be
extended to the magnetohydrodynamic configuration of the VKS experiment
\cite{VKS-PRL-2007}. 
The major component of our algorithm is the poloidal-toroidal 
decomposition \cite{Marques90,Marques93},
which insures incompressibility by construction, at the price of increasing
the order of the governing equations.  When applied to the Navier-Stokes
equation in a finite cylinder, the resulting system has boundary conditions
which are coupled and of high order.  In a companion article
\cite{otherpaper}, we showed that this system could be reduced to the solution
of a set of nested Helmholtz and Poisson problems
with uncoupled Dirichlet boundary conditions, whose solutions could be
superposed via the influence matrix technique.
The purpose of the present article is to describe our method
for solving these elliptic problems using a spectral representation
which exploits azimuthal symmetry of the system and which
is regular at the cylindrical axis, and to demonstrate the validity
of the resulting hydrodynamic code.

More specifically, it was shown by Marques and 
co-workers \cite{Marques90,Marques93} that the Navier-Stokes equations
\begin{subequations}
\label{eq:varr}
\begin{align}
\pdt \bu+(\bu\cdot \bnabla) \bu &= \iR \Delta \bu - \bnabla p \label{v1}\\
\bnabla \cdot \bu &= 0 \label{v2}
\end{align}
\end{subequations}
in a finite cylinder with boundary conditions
\begin{subequations}
\label{eq:cond-u}
\begin{alignat}{3}
\bu&= r\angvel_\pm\et &\qquad\mbox{at}&\quad z=\pmh,\label{eq:cond-uB1}\\
\bu&= 0 &\qquad\mbox{at}&\quad r=1, \label{eq:cond-uB2}
\end{alignat}
\end{subequations}
and with toroidal and poloidal potentials defined by
\begin{equation}
\bu=\curl\left(\psi \ez\right)+\curl\curl\left(\phi\ez\right)
\label{eq:poltor}\end{equation}
are equivalent to the two scalar equations
\begin{subequations}
\label{eq:potMHD_u}
\begin{alignat}{3}
(\pdt-\iR\lap)\laph\psi&= \source_\psi &\equiv &&\ez&\cdot\curl (\bu\cdot\bnabla)\bu 
\label{eq:potMHD_u1}\\
(\pdt-\iR\lap)\lap\laph\phi&=\source_\phi&\equiv
&-&\ez&\cdot\curl\curl(\bu\cdot\bnabla)\bu
\label{eq:potMHD_u2}
\end{alignat}
\end{subequations}
where $\laph\equiv\frac{1}{r}\pdr r \pdr + \frac{1}{r^2}\pd^2_\theta$,
with boundary conditions
\begin{subequations}
\begin{align}
\ir\pdth\psi+\pdrz\phi=\pdr\psi=\laph\phi=\phi=\pdrz\laph\psi-\ir\pdth\lap\laph\phi=0
&\qquad\text{ at }r=1\\
\psi=0&\qquad\text{ at } r=0\\
\laph\psi +\ir\pdr\left(\angvel_\pm r^2\right)= \pdz\laph\phi = \laph\phi=0
&\qquad\text{ at } z=\pmh
\label{eq:bch}
\end{align}
\label{eq:bcs}
\end{subequations} 
Our article \cite{otherpaper} was devoted to showing how the problem
\eqref{eq:potMHD_u}-\eqref{eq:bcs} can in turn be reduced to a sequence of
five nested parabolic and elliptic problems, each with Dirichlet boundary
conditions:
\begin{subequations}
\begin{alignat}{4}
\left(\pdt-\iR\lap\right) f_\psi &= \source_\psi & f_\psi|_{r=1} &= \sigma_f 
& \qquad f_\psi|_{z=\pmh}&=-\frac{1}{r}\pdr (\angvel_\pm r^2)\label{eq:psievol}\\
\laph \psi &= f_\psi & \qquad\pdr\psi|_{r=1} &= 0  & \qquad\psi|_{r=0}&=0
\label{eq:psibcs}\\
\left(\pdt-\iR\lap\right) g &= \source_\phi & g|_{r=1} &= \sigma_g \label{eq:phievol}
& g|_{z=\pmh}&=\sigma^\pm_g\\
\lap f_\phi &= g & f_\phi|_{r=1} &=0 & f_\phi|_{z=\pmh}&=0\\
\laph \phi &= f_\phi & \phi|_{r=1} &= 0
\end{alignat}
\label{eq:big_tableau}
\end{subequations}
for the two potentials $\psi$, $\phi$ and three intermediate
fields $g$, $f_\psi$ and $f_\phi$.
The influence matrix technique \cite{Tuckerman-divfree}, a generalization of the usual separation into
particular and homogeneous solutions, is used to determine boundary values
$\sigma_f$, $\sigma_g$, $\sigma^\pm_g$ such that the boundary conditions
present in \eqref{eq:bcs} but not in \eqref{eq:big_tableau} are satisfied,
\ie~such that
\begin{subequations}
\begin{alignat}{3}
\pdrz f_\psi-\ir\pdth g &=& 0 &\qquad\mbox{at}\quad r=1\\
\ir\pdth\psi+\pdrz\phi&=&0 &\qquad\mbox{at}\quad r=1\\
\pdz f_\phi &=& 0 &\qquad\mbox{at}\quad z=\pmh 
\end{alignat}
\end{subequations}

In this article, we describe the numerical implementation of this algorithm.
We first present the spatial discretization of the fields, using a set of
basis functions \cite{Matsushima} which is regular at the cylindrical axis
$r=0$, and regularizing the discontinuous boundary conditions at the corners
$r=1,z=\pm h/2$.  We then explain the methods we have used for solving equations
\eqref{eq:big_tableau}, in particular for stably and economically
solving the Helmholtz problems resulting from time discretization
and for evaluating the nonlinear terms.
Finally, we describe the validation of the implementation, comparing
results from our code to an analytic polynomial solution and to
previously published two- and three-dimensional test cases.

\section{Spatial and temporal discretization}
\label{sec:spatial_discretization}

The spectral discretization that we use is
\begin{equation}
f(r,\theta,z)=\sum_{m=-\halfM}^{\halfM}f^m(r,z)
e^{im\theta}=\sum_{m=-\halfM}^{\halfM}\,\sum_{k=0}^{K-1}\,
\sum_{
\genfrac{}{}{0pt}{}{n=|m|} {n+m\text{ even}}}
^{2N-1}f_{kn}^m
e^{im\theta}\mQ_n^m(r)\mT_k\left(\frac{2z}{h}\right)
\label{eq:fourdecomp}
\end{equation}
In \eqref{eq:fourdecomp}, we do not introduce new notation for Fourier
coefficients, or for coefficients in the 3D tensor-product basis, using
instead the number and type of superscripts and subscripts to distinguish
between functions in physical space and spectral space coefficients.

The basis functions in the azimuthal and axial directions are 
standard \cite{Orszag71d,Canuto}:
Fourier modes $e^{im\theta}$ and Chebyshev polynomials $\mathcal{T}_k(2z/h)$,
respectively.  
In cylindrical geometries, it is the radial direction which is
most problematic and on which we will focus.  To represent this direction, we
use the basis functions $\mathcal{Q}_n^m(r)$ developed by Matsushima and
Marcus \cite{Matsushima}. In sections \ref{sec:regularity_origin} and
\ref{sec:cornersing}, we will discuss the means by which we impose regularity
at the origin $r=0$ and at the corners $r=1$.

\subsection{Regular basis of radial polynomials}
\label{sec:regularity_origin}

A function $f$ on the disk is analytic at the origin if the radial dependence
of the Fourier coefficient $f^m(r)$ multiplying the Fourier mode $e^{im\theta}$
is of the form
\begin{equation}
f^m(r,z)=\sum^\infty_
{\genfrac{}{}{0pt}{}{n=|m|} {n+m\text{ even}}}
\alpha_n^m(z) r^m
= \alpha_m^m(z) r^m + \alpha^m_{m+2}(z)r^{m+2} + \cdots =r^m p(r^2)
\label{eq:polarsing}\end{equation}
where $p$ is a polynomial.
Examples of functions which violate \eqref{eq:polarsing} are
given in figures \ref{fig:coords_sing}, \ref{fig:coords_sing_order}
and \ref{fig:clustering} and compared with regular functions obeying
\eqref{eq:polarsing} on the right of each figure.

\begin{psfrags}
\psfrag{sing3d}[r][][0.6]{$f(r,\theta)=r$}
\psfrag{reg3d}[r][][0.6]{$f(r,\theta)=r^2$}
\begin{figure}%[h]
	\centering
		\includegraphics[width=0.45\textwidth]{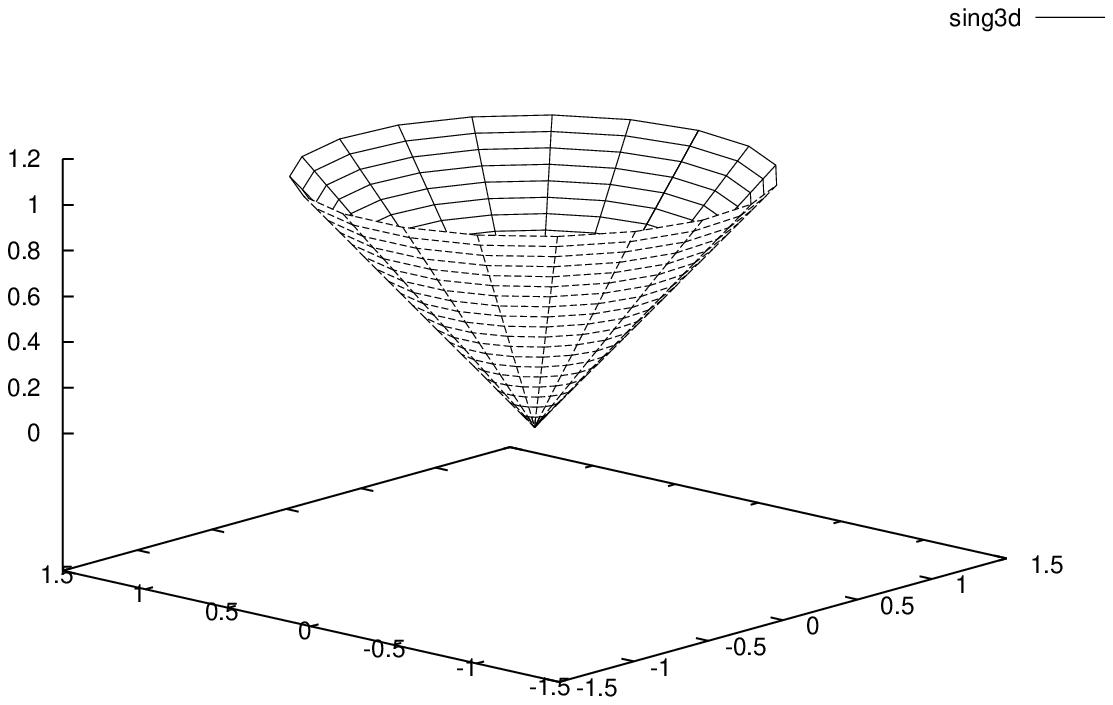}
		\includegraphics[width=0.45\textwidth]{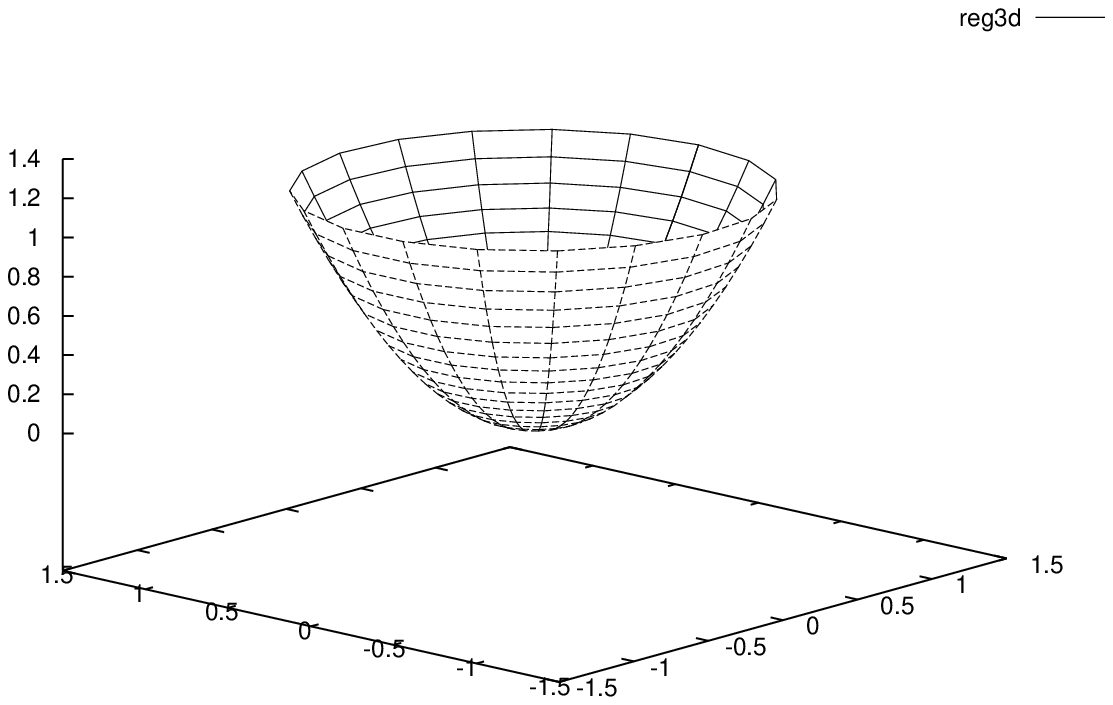}
	\caption{Coordinate singularity effects: parity mismatch. 
Left: $f(r,\theta)=r$. Right: $f(r,\theta)=r^2$.}
	\label{fig:coords_sing}

%\end{figure}
%\end{psfrags}
\vspace*{.5cm}
%\begin{psfrags}
\psfrag{funsig1}[cb][][0.8]{$f(r,\theta)=r^2+\frac{4}{3}\cos(4\theta)$}
\psfrag{funsig2}[cb][][0.8]{$f(r,\theta)=r^2+\frac{4r^2}{3}\cos(4\theta)$}
\psfrag{funreg}[cb][][0.8]{$f(r,\theta)=r^2+\frac{4r^4}{3}\cos(4\theta)$}
%\begin{figure}%[h]
	\centering
		\includegraphics[width=0.32\textwidth]{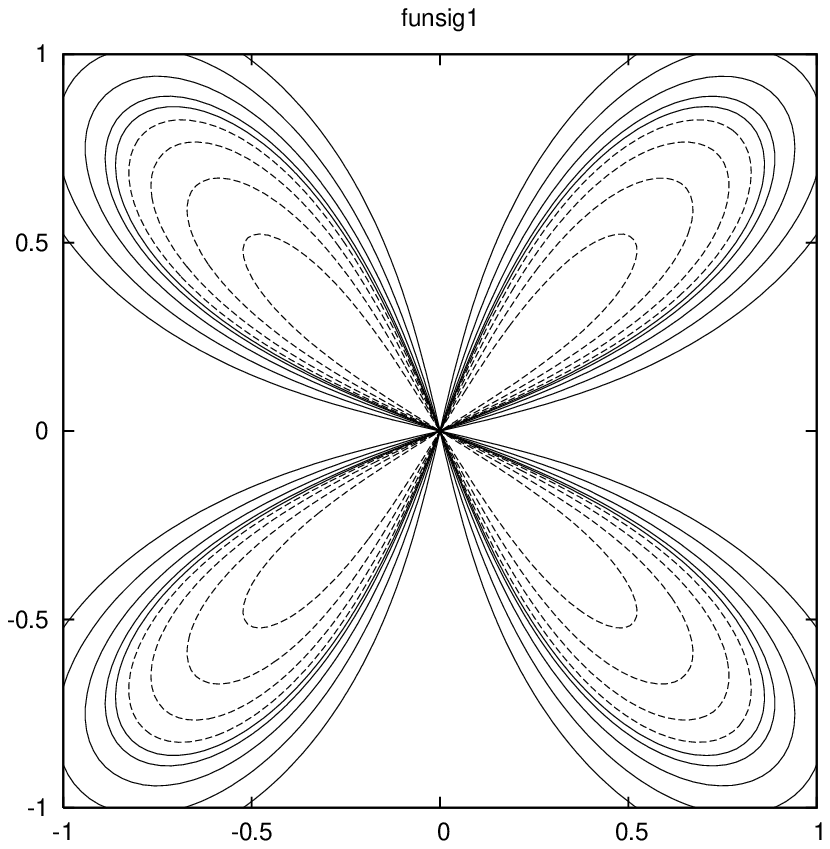}
		\includegraphics[width=0.32\textwidth]{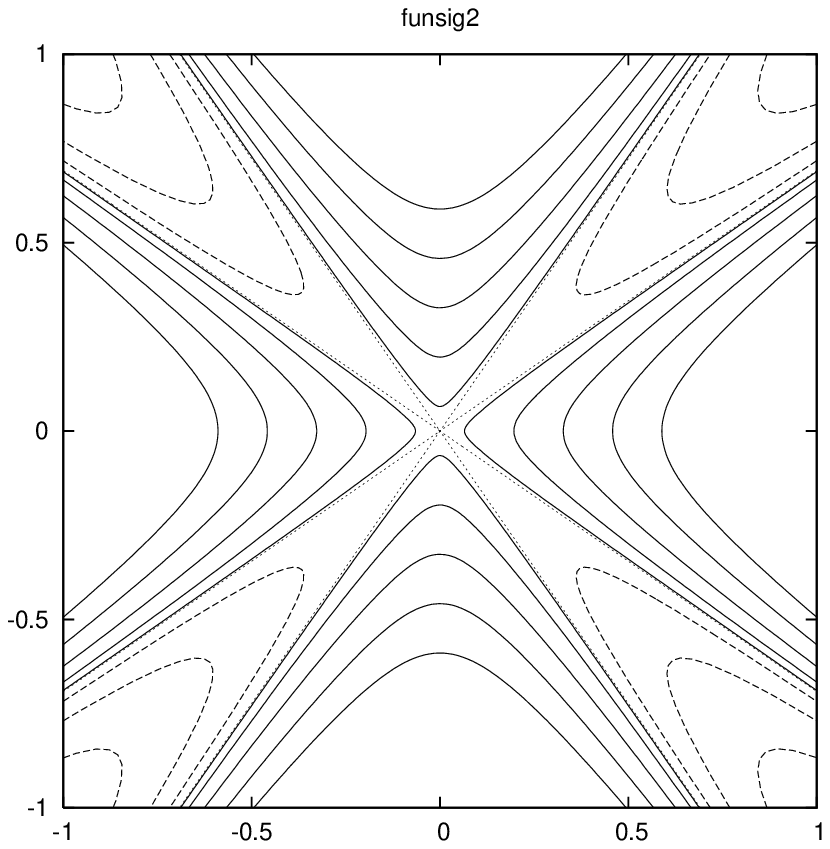}
		\includegraphics[width=0.32\textwidth]{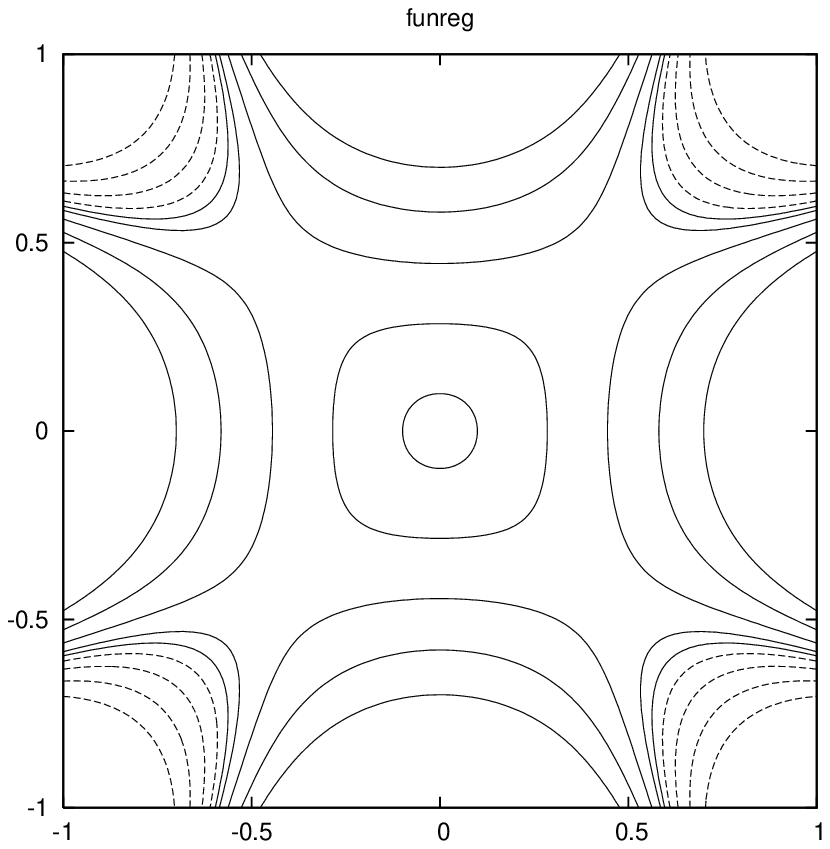}
	\caption{Coordinate singularity effects. Left: discontinuity of
	  value. Middle: discontinuity of Laplacian. Right: regular function.}
	\label{fig:coords_sing_order}
%\end{figure}
%\end{psfrags}
%\begin{psfrags}
\psfrag{g2}[r][r][0.6]{$r^2+32r^2\cos(2\theta)=0.5$}
\psfrag{g4}[r][r][0.6]{$r^2+32r^2\cos(4\theta)=0.5$}
\psfrag{g8}[r][r][0.6]{$r^2+32r^2\cos(8\theta)=0.5$}
\psfrag{g16}[r][r][0.6]{$r^2+32r^2\cos(16\theta)=0.5$}
\psfrag{f2}[r][r][0.6]{$r^2+32r^2\cos(2\theta)=0.5$}
\psfrag{f4}[r][r][0.6]{$r^2+32r^4\cos(4\theta)=0.5$}
\psfrag{f8}[r][r][0.6]{$r^2+32r^8\cos(8\theta)=0.5$}
\psfrag{f16}[r][r][0.6]{$r^2+32r^{16}\cos(16\theta)=0.5$}
%\begin{figure}%[h]
	\centering
		\includegraphics[width=0.45\textwidth]{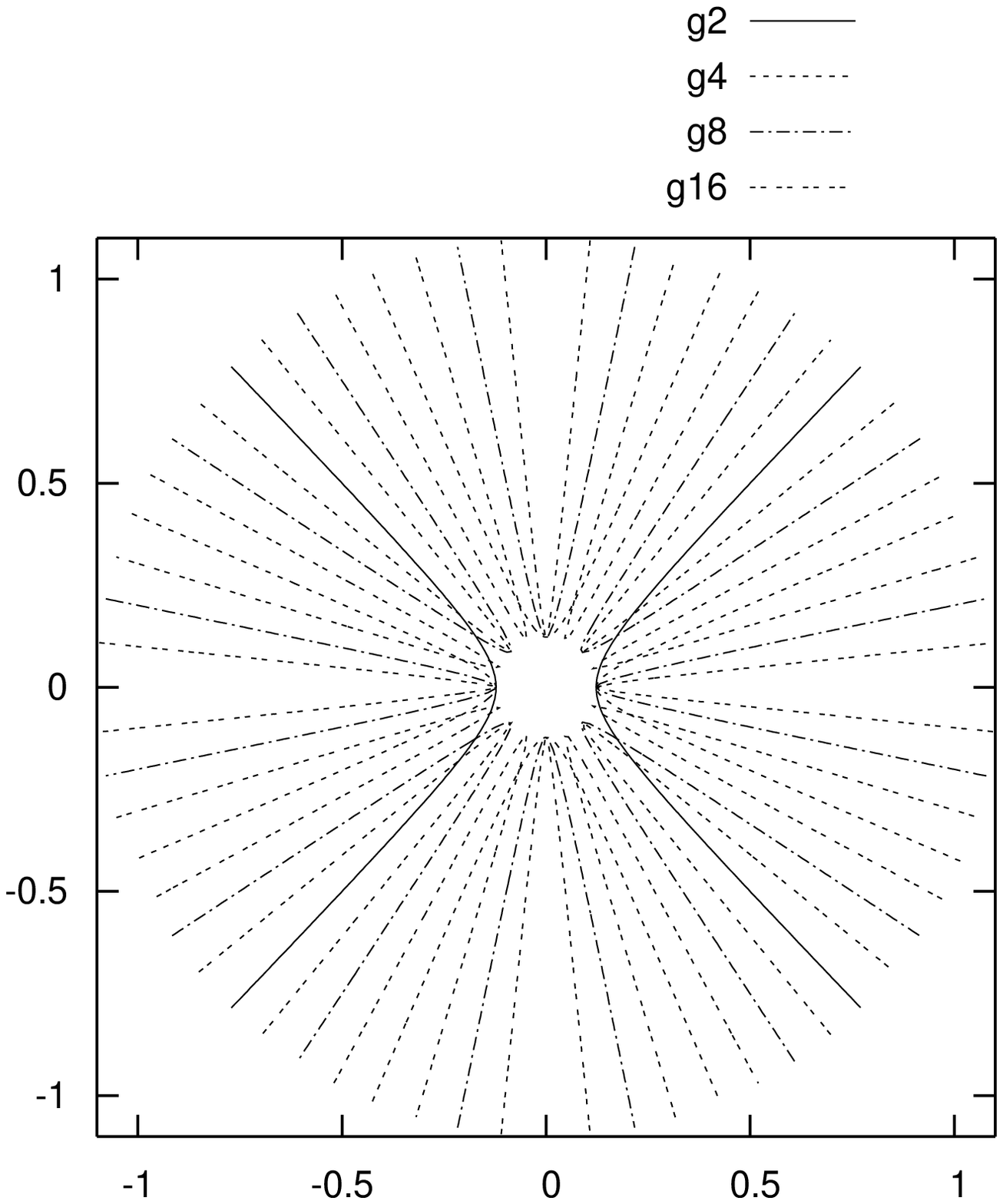}
		\includegraphics[width=0.45\textwidth]{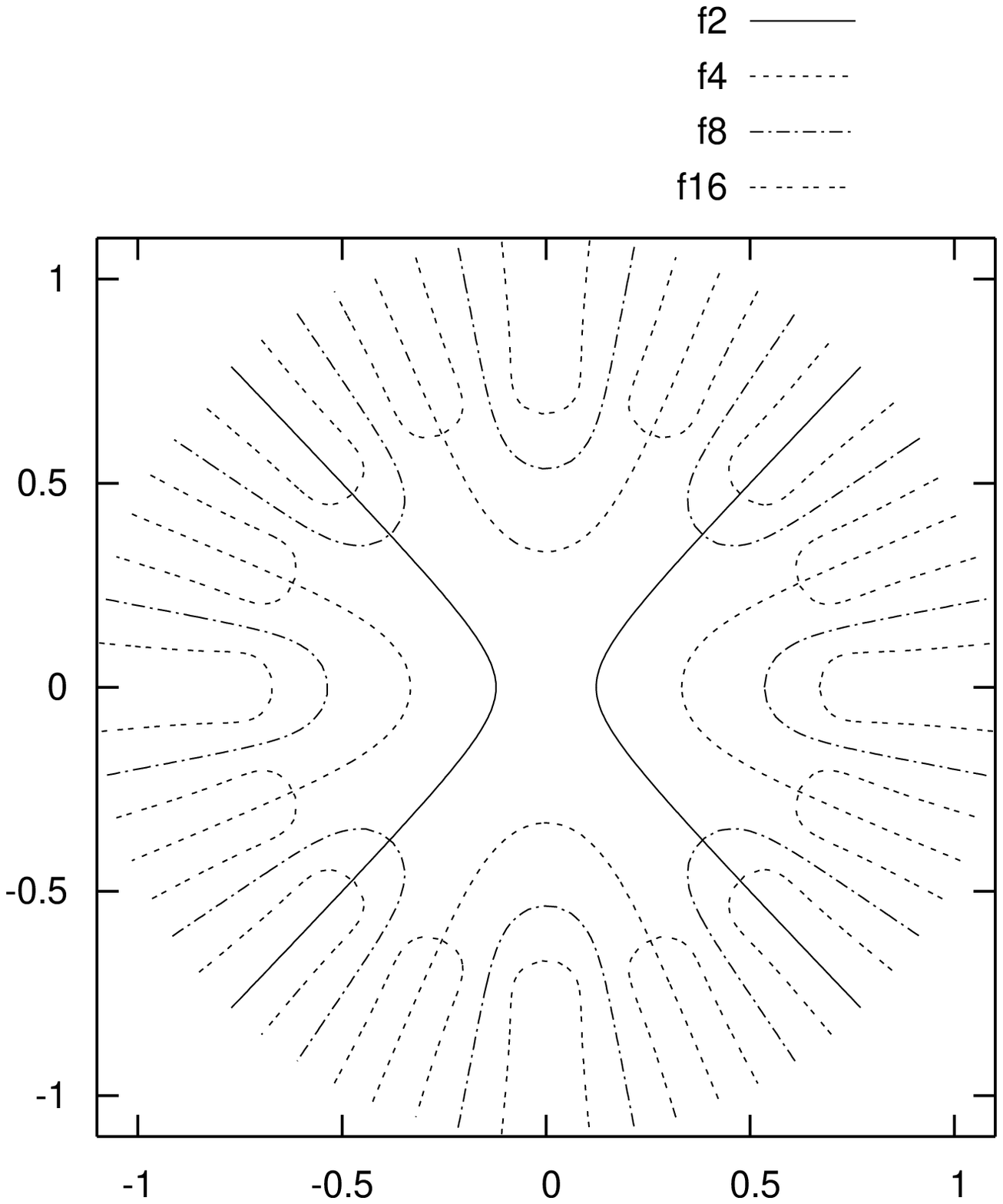}
	\caption{Clustering effect -- contours of $f^m(r,\theta)=0.5$. Left:
	$f^m(r,\theta)=r^2+32r^2\cos(m\theta)$, which is not smooth at 0 for $m>2$.
	Right: $f^m(r,\theta)=r^2+32r^m\cos(m\theta)$, which is smooth at 0.}
	\label{fig:clustering}
\end{figure}
\end{psfrags}

Various approaches used in spectral methods to represent functions 
in polar coordinates are surveyed by Boyd \cite{BoydBook,Boyd99} and
by \citet{Canuto}.
A common practice \cite{LopezMarquesShen2001,Tuckerman-divfree,OrszagPatera,Priymak}
has been to impose some degree of continuity, such as $C^3$, but not 
complete analyticity $C^\infty$.
Although basis functions which are not analytic at the origin
generally do not pollute the fields, retaining such functions wastes
CPU time and memory which could be better used to increase resolution.

The condition \eqref{eq:polarsing} is stated in terms of monomials, the use of
which is excluded because of their poor numerical properties.  The polynomial
basis developed by \citet{Matsushima} respects these conditions and yet is
numerically well-conditioned. These polynomials
$\mathcal{Q}_n^m(\alpha,\beta;r)$ are solutions to the singular
Sturm-Liouville equation:
\begin{equation}
\left(\frac{\left(1-r^2\right)^{1-\alpha}}{r^\beta}\frac{d}{dr}\left(\left(1-r^2\right)^\alpha
r^\beta\frac{d}{dr}\right)-\frac{|m|(|m|+\beta-1)}{r^2}+n(n+2\alpha+\beta-1)\right)\mathcal{Q}_n^m(\alpha,\beta;r)=0,
\label{eq:sturm}\end{equation}
defined over $r\in[0,1]$.
In \eqref{eq:sturm}, $0\leq |m| \leq n$, $\alpha\in[0,1]$
and $\beta$ is a positive integer.
With the special choice $\alpha=\beta=1$,
$\mathcal{Q}^m_n(1,1;r)$ are related to Legendre and shifted Jacobi
polynomials used by \citet{Leonard}; similar functions were also derived
by \citet{Verkley}.
The functions $\mathcal{Q}^m_n(\alpha,\beta;r)$ are complete and orthogonal
over $[0,1]$ with respect to the inner product:
\begin{equation}
\int_0^1\mathcal{Q}_n^m(\alpha,\beta;r)\mathcal{Q}_{n'}^m(\alpha,\beta;r)
\frac{r^\beta}{\left(1-r^2\right)^{1-\alpha}} =
I_n^m(\alpha,\beta)\delta_{nn'}
\end{equation}
The $n^{\rm th}$ order polynomials $\mathcal{Q}_n^m(\alpha,\beta;r)$
have the following explicit expression:
\begin{equation}
\mathcal{Q}_n^m(\alpha,\beta;r)\equiv
\sum_{p=0}^{\frac{n-|m|}{2}}\frac{(-1)^{p+\frac{n-|m|}{2}}
\Gamma\left(\frac{n+|m|+\gamma-1}{2}+p\right)
\Gamma\left(\frac{2|m|+\beta+1}{2}\right)}{\Gamma(p+1)
\Gamma\left(\frac{n-|m|}{2}-p+1 \right)
\Gamma\left(\frac{2|m|+\beta+1}{2}+p\right)
\Gamma\left(\frac{2|m|+\gamma-1}{2}\right)}
r^{|m|+2p},
\label{eq:radial_poly_def}
\end{equation}
but they, as well as the normalizing coefficients $I_n^m(\alpha,\beta)$, can
be calculated in $O(n-|m|)$ operations using recursion relations given by
\citet{Matsushima}.  Recursion relations also exist for the operators
\begin{equation}
\left\{r\pdr, r^2, (r\pdr)^2-m^2,(r\pdr)^2+\lambda r^2\right\}
\label{eq:parity_ops}
\end{equation}
expressed in the $\mathcal{Q}_n^m$ polynomial basis,
meaning that for any of the operators $H$ in \eqref{eq:parity_ops},
there exist banded matrices $L$ and $R$ such that $H = R^{-1} L$.
Thus,
\begin{equation}
H f = g \Longleftrightarrow L f = R g \label{eq:rec3}
\end{equation}
reducing the time for multiplication by $H$ or $H^{-1}$ from
quadratic to linear in the number of radial modes or gridpoints.  The
existence of recursion relations is a general property of differential
operators represented in polynomial bases, for reasons explained by
\citet{Tuckerman-banded}.
Recursion relations will be further discussed in section \ref{sec:viscous}.

The radial function $f^m(r)$ associated with Fourier mode $m$ and its
coefficients $f^m_n$ in the polynomial basis are related by the
transform pair:
\begin{subequations}
\begin{align}
f^m(r)&=\sum_{\genfrac{}{}{0pt}{}{n=|m|} {n+m\text{ even}}}
^\infty f_n^m \mathcal{Q}_n^m(r)
\approx\sum_{\genfrac{}{}{0pt}{}{n=|m|} {n+m\text{ even}}}
^{\hat{N}}f_n^m\mathcal{Q}_n^m(r)\\
f_n^m &= \int_0^1 dr \: w(r) f^m(r)\mathcal{Q}_n^m(r)/I_n^m
\approx \sum_{i=0}^{N-1} w_i f^m(r_i)\mathcal{Q}_n^m(r_i)/I_n^m
\end{align}
\label{eq:radial_decomposition}
\end{subequations}
In \eqref{eq:radial_decomposition}, the order of the polynomial
expression is $\hat{N}\equiv 2N-2$ or $\hat{N}\equiv 2N-1$ according to
whether $m$ is even or odd.  The collocation points $\{r_i\}$ for
Gaussian quadrature are computed
numerically as the roots of the first neglected $m=0$ polynomial
$\mathcal{Q}^0_{\hat{N}+2}$ if the boundary points are to be excluded;
otherwise, they are the roots of a slightly more complicated expression
\cite{Matsushima}. Once the $\{r_i\}$ are determined,
the weights $\{w_i\}$ are computed by recursion relations \cite{Matsushima}.

Equations \eqref{eq:radial_decomposition} specify $N_m\equiv N-[m/2]$ 
coefficients from
values at $N$ quadrature points via a rectangular matrix.  Since the basis is
orthonormal, the inverse transformation is obtained from the transpose of this
rectangular matrix.  The smaller size of the spectral representation compared
to the grid representation is a consequence of the fact that the functions in
\eqref{eq:radial_decomposition} are not arbitrary functions of $r$, but belong
to the restricted subspace of functions obeying the regularity conditions
\eqref{eq:polarsing}.

\subsection{Regularization of the corners}
\label{sec:cornersing}

The boundary conditions on $u_\theta$ stated in \eqref{eq:cond-u} are
\begin{subequations}
\begin{alignat}{3}
\quad u_{\theta}(r,\theta)&= r\angvel_{\pm} &\qquad\mbox{at}&\quad z=\pmh\label{eq:bc_disks}\\
u_{\theta}(\theta,z)&=0 &\qquad\mbox{at}&\quad r=1 \label{eq:bc_cyl}
\end{alignat}
\label{eq:bc_uth}
\end{subequations}
The equations have been nondimensionalized using the radius as the 
unit of length and the inverse angular velocity 
$1/\angvel_-$ (with $\angvel_- > 0$ and $\vert\angvel+\vert\leq
\angvel_-$) as the unit of time.
Therefore, $\angvel_-=1$ and $-1 \leq s=\angvel_+ \leq 1$,
in particular $\angvel_+$ is $0$, $+1$ or $-1$, for the
rotor-stator, exactly corotating, or exactly counter-rotating configurations, 
respectively. (It is also possible to set the velocity on the cylinder
to some other constant value instead of 0; for example
to simulate the flow in a rotating cylinder,
we would set $u_{\theta}|_{r=1}(\theta,z)=\omega_+=\omega_-=1$.
Here, we focus instead on the difficulties in implementing 
the discontinuous boundary conditions \eqref{eq:bc_uth}.)

Boundary conditions \eqref{eq:bc_uth} are discontinuous at the corner points
$r=1,\ z=-h/2$ and possibly $r=1,\ z=+h/2$: one or both disks rotate while the
lateral boundary remains fixed.  Mathematically, a PDE with a finite number of
singular points can have a solution which is smooth except at these points.
However, spectral methods then do not converge exponentially because series of
smooth functions cannot converge uniformly to a discontinuous solution. If
nothing is done to prevent it, the Gibbs phenomenon will lead to spurious
oscillations which propagate into the whole domain from the neighborhood of
the singularity. For finite difference methods, the discontinuity will affect
only the neighborhood of the singular point, on the order of the grid
interval, and therefore does not pose a severe problem.  Finite volume methods
have a local integral formulation and so the discontinuity presents an even
less serious problem.  The filtering intrinsic to local methods is, however,
intrinsically related to the high numerical diffusion which in turn makes
local methods less precise.

In some cases, even local methods do not sufficiently filter singularities.
\citet{Georgiou} discuss the issue of spurious oscillations
in the context of finite element methods.
In the solutocapillary problem studied by \citet{Witkowski-Walker},
the authors were required to explicitly filter the solution to achieve
acceptable convergence even in a finite difference calculation.

\begin{figure}[h]
\centering
%\begin{minipage}[c]{0.45\textwidth}
a)\includegraphics[width=0.25\textwidth]{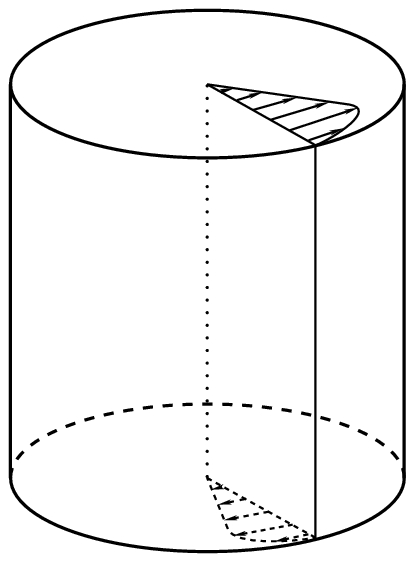}
%\end{minipage}
\qquad\qquad\qquad
%\begin{minipage}[c]{0.45\textwidth}
b)\includegraphics[width=0.25\textwidth]{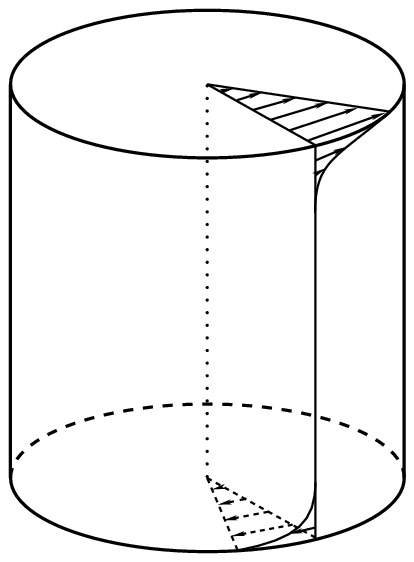}
%\end{minipage}
\caption{Regularized profiles used for elimination of the discontinuous
boundary conditions at the cylinder corners. a) regularization on upper and
bottom disks using \eqref{eq:reg}.  b) regularization on
lateral bounding cylinder using \eqref{eq:reg_prof_uth_lat}.}
\label{fig:reg_prof}
\end{figure}

Spectral methods must always explicitly filter strong singularities
like \eqref{eq:bc_disks}-\eqref{eq:bc_cyl}.
We have chosen to do this by
approximating the discontinuous function at the boundary by a steep but smooth
profile.  This procedure can be justified by arguing that we are not
interested in finding the solution to the singular problem.  In any real
experiment, the boundary conditions are not discontinuous: a small gap must
necessarily exist between the rotating disks and the stationary lateral
boundary.  In our algorithm, we replace the constant angular velocities 
in \eqref{eq:bc_disks} by continuous functions,
as illustrated by figure \ref{fig:reg_prof}a.  Two possibilities are:
\begin{subequations}
\begin{alignat}{3}
u_\theta(r,\theta) &=r\left(1-e^{\frac{r-1}{\delta}}\right)\angvel_\pm
&\qquad\mbox{at}&\quad z=\pmh \label{eq:reg_exp}\\
u_\theta(r,\theta)&=r(1-r^\mu)\angvel_\pm &\qquad\mbox{at}&\quad z=\pmh
\label{eq:reg_poly}
\end{alignat}
\label{eq:reg}
\end{subequations}\\
where $\mu$ is an arbitrary but large even integer (e.g. $\mu=10$),
and $\delta$, $\angvel_\pm$ are constants.
In order to be represented in a radial polynomial basis, the exponential
regularization \eqref{eq:reg_exp} must be approximated by a polynomial, so
that both expressions above are effectively polynomials.
The steepness of the profiles are adjusted by varying
$\delta$ or $\mu$, whose possible values are limited by the radial
polynomial order $N$.

The advantage of exponential regularization is that, for a given $N$,
a steeper profile can be achieved by the polynomial approximation to
\eqref{eq:reg_exp} than by the polynomial \eqref{eq:reg_poly},
as can be seen in figure \ref{fig:reg_prof_r}a.
In this way, the deviation from the idealized profile is minimized,
without unduly increasing the radial resolution $N$ over
that required to resolve the field in the interior.
The polynomial approximation to \eqref{eq:reg_exp}
differs from 1 by more than 10\% only over the small range $1-2\delta\
{\scriptstyle\lessapprox}\,r\,<1$.

Regularization of the boundary condition imposes a lower bound on the spectral
resolution -- the spectral approximation must be able to represent the
regularization profiles smoothly. In fact, \citet{LopezShen1998} observed that
the actual resolution should be approximately twice the minimal resolution
sufficient for representing the regularization profiles, due to generation of
higher wavenumbers by the nonlinear terms. As was shown by \citet{LopezShen1998}, for comparing
with results obtained using different methods and for benchmark purposes it is
sufficient to use $\delta\approx 0.005$.  In practice, we use
$0.005<\delta<0.05$, illustrated in figure \ref{fig:reg_prof_r}b.

Care must be taken in approximating \eqref{eq:reg_exp} by a
polynomial expression, in order to satisfy all of the conditions that
we require of $u_\theta$ and of $\psi$.  The procedure we use is as
follows.  We evaluate \eqref{eq:reg_exp} on the collocation points.  Since
$u_\theta$ is an odd function of $r$ (see \eqref{eq:vector}), 
we apply the transform
\eqref{eq:radial_decomposition} using the (odd) polynomials associated with
the $m=1$ Fourier mode.  The basis of odd radial functions insures that this
approximation to $u_\theta$ on the bounding disks is zero at $r=0$, while the
use of Gauss-Radau collocation points which include the cylinder boundary
insures that that it is also zero at $r=1$.  This approximation to
$u_\theta$ is integrated over $r$ to obtain an approximation to $\psi$, with
the integration constant chosen in order to satisfy $\psi=0$ at $r=0$. 

\begin{psfrags}
\psfrag{regfun1}[r]{$\delta=0.05$}
\psfrag{regfun2}[r]{$\delta=0.005$}
\begin{figure}[hc]
\centering
a)\includegraphics[width=7cm]{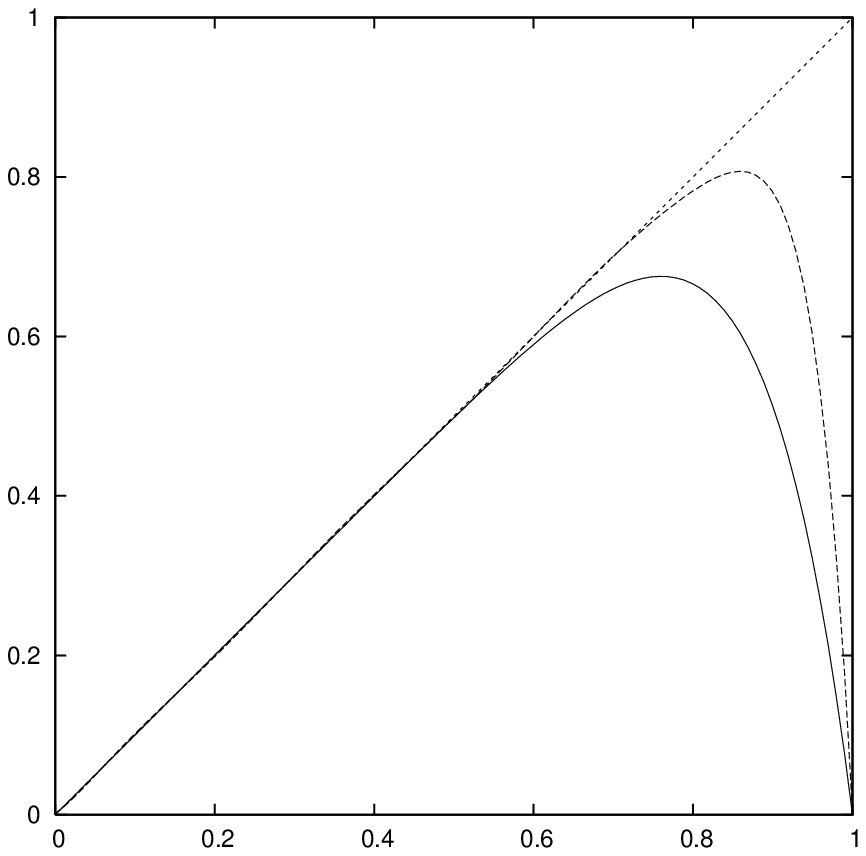}
b)\includegraphics[width=7cm]{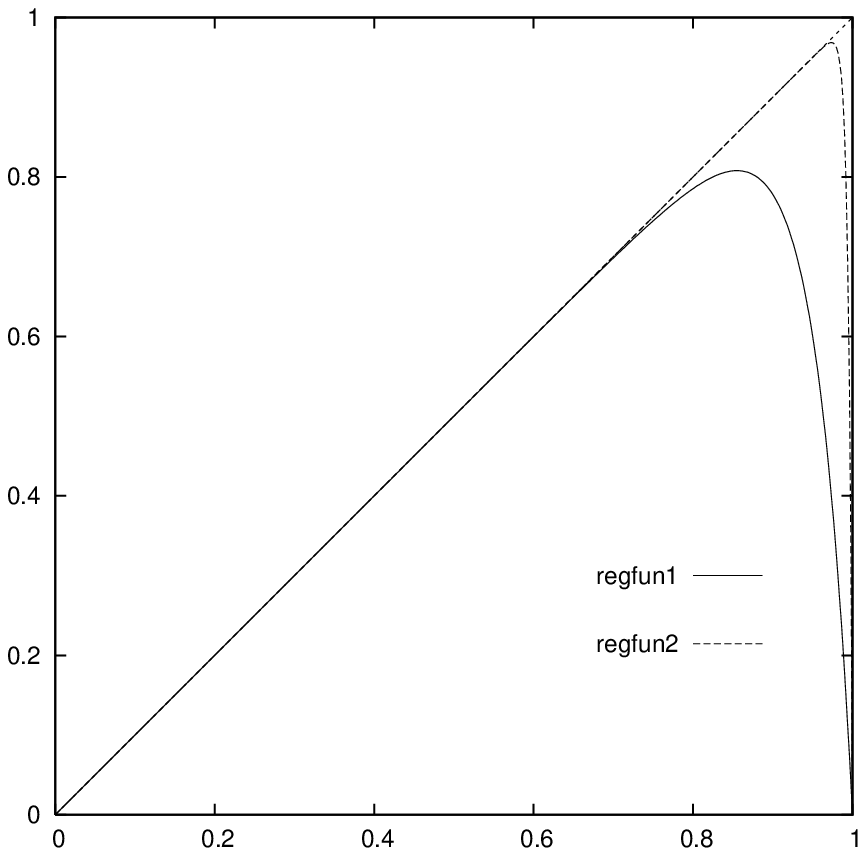}
\caption{
a) Comparison between regular profiles represented by polynomials of
order $r^9$. Solid curve: $u_{\theta}|_{z=\pmh}(r)= r\left(1-r^8\right)$.
Dotted curve: $u_{\theta}|_{z=\pmh}(r)\approx
r\left(1-e^{\frac{r-1}{\delta}}\right)$. Inset illustrates result 
b) Comparison between exponential profiles 
$\left.u_{\theta}\right|_{z=\pmh}(r)= r\left(1-e^{\frac{r-1}{\delta}}\right)$
for different $\delta$.}
\label{fig:reg_prof_r}
\end{figure}
\end{psfrags}

Another choice is to apply a filter to the lateral boundary,
replacing \eqref{eq:bc_cyl} by
\begin{equation}
u_{\theta}(z) = 
\angvel_+e^{-\left(1-\frac{2z}{h}\right)/\delta} + 
\angvel_-e^{-\left(1+\frac{2z}{h}\right)/\delta}
\label{eq:reg_prof_uth_lat}\qquad\mbox{at}\quad r=1
\end{equation}
while keeping the uniform angular velocity profiles \eqref{eq:bc_disks}
on the disks unchanged; see figure \ref{fig:reg_prof}b.
This kind of regularization is similar to that of \citet{LopezShen1998}. 

Other, quite different, approaches to the treatment of singularities exist.
One is singularity subtraction.  The form of the singular part of the solution
is determined analytically, and the solution is written as a sum of the
singular solution and an unknown regular part. Only the regular part is
treated numerically. The effect of the singular solution on the numerical one
can be filtered down to the scales representable by the spatial
resolution. The main advantage of this method is that it recovers the
convergence of the scheme and at the same time approaches the exact
solution. Recent applications of this method are to the driven cavity problem
\cite{Schultz} and to injection of fluid into a cylindrical channel
\cite{Botella}. The results obtained are generally of high precision and often
provide a benchmark for a particular problem.  The main drawback is that it
requires knowledge of the solution near the singular point. For the 2D driven
cavity problem, the nature of the singularity was given by \citet{Dean} and
\citet{Moffatt} for a Stokes flow.  For most inertial (Navier-Stokes) flows,
and for 3D flows \citet{Hills}, the analytic form of the 
singular solution is unknown.  
Note that even when the velocity boundary conditions 
are continuous, lower-order singularities of purely geometric
origin are present at the corners.
We will return to this in section \ref{sec:non-poly}.

Another approach is to derive a physically justified model which is no longer
singular.  Methods from molecular dynamics reflect the microscopic nature of
the fluid at the smallest scales but are very hard to adapt to problems
containing both large and small scales.  Several continuous (macroscopic)
approaches have been proposed as a compromise between a continuous and a
molecular description. These all introduce a spatially limited physical effect
which effectively removes the singularity. In this category are methods based
on variable slippage, as well as the surface viscosity or dynamic surface
tension applicable to free-surface problems. A comprehensive review of
physically justified models, as well as other regularization techniques, is
provided by Nguyen and co-workers \cite{SebBook,SebThese}.

\subsection{Time integration}
\label{sec:time}

We use an implicit scheme for the linear diffusive terms while treating the
nonlinear terms explicitly.  Spectral methods require that the coefficients
representing the solution decay with their index or wavenumber. The nonlinear
term in the Navier-Stokes equation can be seen as a generator and amplifier of
high wavenumbers while the viscous term damps these high wavenumbers. The
intensity of this damping depends on the particular time-integration scheme
and on the way the Laplacian is evaluated and must be strong enough to oppose
the effect of the nonlinear term.  In our case, high-wavenumber modes are
needed to represent both the thin boundary layers created near the rotating
disks and the steep regularized boundary profile.  Fortunately, these effects
are most pronounced in the proximity of the boundaries, where the axial
Chebyshev and radial polynomial grid is finest.  However, in the
counter-rotating case, the central shear layer can also require high
wavenumber modes in order to be well-represented.

We use the first-order backward Euler scheme for linear terms because it
attenuates high frequencies faster than all other methods.  Tests
performed with the Crank-Nicolson method confirmed that for this scheme,
nonlinear simulation was unstable even for quite small Reynolds numbers
$Re\approx300$. This behavior was also observed for the von K\'arm\'an flow by
\citet{Speetjens} and \citet{LopezMarquesShen2001} and by
\citet{Marcus84} in Taylor-Couette flow.
The nonlinear term is treated by the second-order explicit Adams-Bashforth
scheme, so that the backward Euler/Adams-Bashforth time-integration scheme for
the potentials $\psi$ and $\phi$ takes the following form:
\begin{subequations}
\begin{align}
(I-\dt\iR\lap)\lap_h\psi^{n+1} &= \lap_h\psi^{n} + \frac{\Delta
    t}{2}\left(3 S^n_{\psi}-S^{n-1}_{\psi}\right) \equiv rhs_\psi^n
\label{eq:timestep_psi}\\
(I-\dt\iR\lap)\lap\lap_h\phi^{n+1} &= \lap\lap_h\phi^{n} +
    \frac{\Delta t}{2}\left(3 S^n_{\phi}-S^{n-1}_{\phi}\right) \equiv rhs_\phi^n
\label{eq:timestep_phi}
\end{align}
\label{eq:timestep}
\end{subequations}
where $\source_\psi$ and $\source_psi$ are defined in \eqref{eq:potMHD_u}.
Equations \eqref{eq:timestep}
can be written as the nested system of equations:
\begin{subequations}
\begin{alignat}{4}
%\left(I-\dt\iR\lap\right) f_\psi &= rhs_\psi^n & \qquad\laph \psi^{n+1} &= f_\psi  \\
%\left(I-\dt\iR\lap\right) g &= rhs_\phi^n & \qquad\lap f_\phi &= g & \qquad\laph \phi^{n+1} &= f_\phi
\left(I-\dt\iR\lap\right) f_\psi &= rhs_\psi \label{eq:system1}\\
\laph \psi &= f_\psi  \label{eq:system2}\\
\left(I-\dt\iR\lap\right) g &= rhs_\phi \label{eq:system3}\\
\lap f_\phi &= g \label{eq:system4}\\
\laph \phi &= f_\phi \label{eq:system5}
\end{alignat}
\label{eq:big_tableau_bis}
\end{subequations}%
As explained in section \ref{sec:intro}, the boundary conditions
imposed on \eqref{eq:big_tableau_bis} are Dirichlet conditions
with boundary values calculated via the influence matrix in 
such a way as to satisfy the more complicated coupled boundary conditions
given in \eqref{eq:bcs}.

The maximal time step $\dt$ depends on the Reynolds number. Typically starting
from state $\bu=0$ requires a $\dt$ which is 4--10 times smaller than that
which can be used for evolving a fully developed state at the same Reynolds
number. This is because the state $\bu=0$ is incompatible with the boundary
condition \eqref{eq:bc_disks}. In the first few iterations a boundary layer is
created near the rotating cylinder lids, requiring higher spatial
resolution. This can be avoided by performing about 100 initial steps of the
linear Stokes solver, \ie~ with $\source_\psi=\source_\phi=0$.

\begin{table}[!h]
\centering
\begin{tabular}{cccc}
\toprule
$Re$ & configuration & $\dt$ & resolution ($M\times K\times N$) \\
\midrule
< O(500) & 2D & 0.05 - 0.1 & $1\times 32\times 16$ \\
500 - 1000 & 2D & 0.02 - 0.05 & $1\times 64\times 32$ \\
1000 - 3000 & 2D & 0.01 - 0.02 & $1\times 96\times 48$ \\
3000 - 5000 & 2D & 0.005 - 0.01& $1\times 128\times 64$\\
5000 - 10000& 2D & 0.001 - 0.0025& $1\times 180\times 90$\\\hline
< O(500) & 3D & 0.04 - 0.1 & $8\times 64\times 32$ \\
500 - 1000 & 3D & 0.01 - 0.04 & $16\times 80\times 40$ \\
1000 - 3000 & 3D & 0.025 - 0.01 & $32\times 100\times 60$ \\
3000 - 5000 & 3D & 0.001 - 0.0025& $(64$-$96)\times 128\times 80$\\
\bottomrule
\end{tabular}
\caption{Typical values of timestep $\dt$ and azimuthal ($M$), axial
($K$) and radial ($N$) resolution, for different
configurations. Typical steepness of the regularization profile is
$\delta\approx 0.01$}
\label{tab:timesteps}
\end{table}

In Table \ref{tab:timesteps} we present the values of $\dt$ and spatial
resolutions typically used for performing nonlinear simulations for different
values of $Re$.

\subsection{Viscous terms\label{sec:viscous}}

We now describe the way in which the Helmholtz and Poisson problems with
Dirichlet boundary conditions in \eqref{eq:big_tableau} or
\eqref{eq:big_tableau_bis} are solved.  The azimuthal Fourier representation
in \eqref{eq:fourdecomp} decomposes each 3D elliptic problem into a set of 2D
problems, each of which is associated to a single azimuthal Fourier mode $m$.
The reflection symmetry in $z$ leads to further decoupling between the set of
modes that are symmetric or antisymmetric in $z$.  Each of the $2M$ resulting
elliptic problems corresponds to a single azimuthal Fourier mode $m$ and 
axial parity $p\in\{s,a\}$, within which each equation corresponds to a value of
$(k,n)$, the indices of the axial and radial basis functions.  The number of
axial modes of each parity is $K/2$, and the number the number of radial basis
functions corresponding to Fourier mode $m$ is $N_m\equiv N-[\frac{m}{2}]$; in
the remainder of this section, we will take $m=0$, so that $N_m=N$.

In the full cylinder, regularity at the axis serves as one of the boundary
conditions and is imposed by the use of the $\mQ_n$ radial basis in
\eqref{eq:fourdecomp}, leaving only the boundary condition at $r=1$ to be
imposed.  In the axial direction, the boundary conditions at the two disks can
be recombined to yield one condition for each parity.  Thus, one radial and
one axial boundary condition remains to be imposed on each 2D problem.  These
are imposed via the $\tau$ method \cite{Orszag71d,Canuto}, so that the
equations corresponding to the highest-wavenumber modes in each direction are
replaced by the boundary conditions.

Gaussian elimination is especially economical for systems resulting from the
spectral discretization of differential equations, whose structure is barely
altered when boundary conditions are imposed via the $\tau$ method, since
recursion relations reduce the solution time from quadratic to linear 
in the number of modes if the resulting systems are diagonally dominant.
This can only be done in one direction however. In geometries with
more than one non-periodic direction, the remaining directions
must be treated by diagonalization.
Here, we treat the axial direction by incorporating the boundary condition via
Schur decomposition into the matrices representing $\pd^2_z$ for each parity
and diagonalizing \cite{Canuto,Haidvogel}, leading to decoupled problems
for each axial eigenvalue $\lambda_z$. 
The operation count at each timestep, dominated by multiplication
by the eigenvector matrix, is quadratic in $K/2$ and linear in $N$.

Thus, the 2D and 3D problems in \eqref{eq:big_tableau_bis}
%\begin{align}
%\laph f &= g \label{eq:laph}\\
%\lap f &= g \label{eq:lap}\\
%(I-\dt\iR \lap) f &= g \label{eq:helm}
%\end{align}
are all decomposed into a set of one-dimensional radial problems:
\begin{align}
H f &\equiv\left[\frac{1}{r}\pdr r\pdr -\frac{m^2}{r^2}+\lambda \right]f = g\\
\intertext{which we will write in practice as}
r^2 H f &\equiv\left[r\pdr r\pdr -m^2+\lambda r^2\right]f = r^2 g
\label{eq:Hprob}\end{align}
where $m$ is the Fourier mode. The scalar $\lambda$ is 0 or $\lambda_z$
for the Poisson problems
%\eqref{eq:laph} or \eqref{eq:lap}, or
\eqref{eq:system2}, \eqref{eq:system4} or \eqref{eq:system5}, or
$\lambda_z+Re/\dt$ for the 
%Helmholtz problem \eqref{eq:helm} 
Helmholtz problems \eqref{eq:system1} and \eqref{eq:system3} 
(with the multiplicative factor $-Re/\dt$ incorporated into $g$).

With $f$ and $g$ represented in the polynomial radial basis
\eqref{eq:radial_decomposition}, a recursion relation exists for
$r^2 H$, as stated in section \ref{sec:regularity_origin}, 
\ie~$r^2 H = R^{-1}L$ with $R$, $L$ banded matrices.
Thus each Helmholtz problem \eqref{eq:Hprob} can be replaced by
\begin{align}
L f = R r^2 g &\equiv Q g\label{eq:rec4}
\end{align}
For nonzero $\lambda$, $L$ is pentadiagonal and $R$ is tridiagonal.  Two
obstacles must be surmounted before \eqref{eq:rec4} can be solved in a time
which is linear in the number of radial modes.  First, $L$ is not
diagonally dominant, so that stable Gaussian elimination would require
pivoting, destroying the banded structure.  Second the radial boundary
condition must be included.  A method which overcomes these difficulties was
presented by \citet{Matsushima}.  Here, we will recast this method in terms of
the Sherman-Morrison-Woodbury formula, which can be shown \cite{Tuckerman-divfree}
to underly a large class of transformations between coupled and uncoupled
systems.

System \eqref{eq:rec4} must be altered in two ways.
First, since the largest element of $L$ is located on the first super-digonal,
permuting its rows leads to a matrix $PL$ which is diagonally dominant, 
but is no longer banded. 
Second, the boundary conditions must be imposed.
The tau method replaces \eqref{eq:Hprob} by
\begin{subequations}
\begin{align}
H  f &= g + \hat{e}_N \tau \label{eq:tau}\\
B^T f &= \beta \label{eq:bc}
\end{align}
\end{subequations}
which effectively discards $g_N$ by adjusting it with the extra unknown
$\tau$, introduced along with the boundary condition \eqref{eq:bc}.
$B^T$ is the row vector which represents the discretized boundary condition 
and $\hat{e}_N$ is the unit vector which selects the $N^{th}$ component.
Multiplying by $r^2$, using \eqref{eq:rec4} and permuting rows leads to:
\begin{align}
r^2 H f = R^{-1} L f & = r^2 g + r^2 \hat{e}_N \tau \nonumber\\
L f & = R r^2 g + R r^2 \hat{e}_N \tau \nonumber\\
P L f & = P Q g + P Q \hat{e}_N \tau 
\end{align}
This is rewritten in matrix form as
\begin{equation}
\left(
\begin{MAT}(@)[2pt]{c2c}
L_{N,i} & -Q_{N,N}\\2
\begin{MAT}(@)[1pt,3.5cm,1.5cm]{c}
L_{lo,i}\\
\end{MAT} & -Q_{lo,N} \\2
B^T & 0 \\
\end{MAT}
\right)
\left(
\begin{MAT}(@)[2pt]{c}
\begin{MAT}(@)[1pt,0.5cm,2.1cm]{c}
f \\
\end{MAT}\\2
\begin{MAT}(@)[1pt,0.5cm,0.4cm]{c}
\tau \\
\end{MAT}\\
\end{MAT}
\right)
=
\left(
\begin{MAT}(@)[2pt]{c}
(Qg)_N \\2
\begin{MAT}(@)[1pt,0.5cm,1.5cm]{c}
(Qg)_{lo} \\
\end{MAT}\\2
\beta\\
\end{MAT}
\right)
\label{eq:badmatrix}
\end{equation}
where the subscript $lo$ refers to all indices lower than $N$.
System \eqref{eq:badmatrix} is solved by using
the Sherman-Morrison-Woodbury formula
\begin{equation}
(A + v w^T)^{-1} = A^{-1} - A^{-1} v (I + w^T A^{-1} v)^{-1} w^T A^{-1}
\label{eq:smw}
\end{equation}
which relates the inverses of two matrices differing by a low-rank
transformation $vw^T$, in particular differing by a few rows or columns.
We define $A$, $v$ and $w^T$ for \eqref{eq:badmatrix} as follows.
The matrix in \eqref{eq:badmatrix} is replaced by another matrix $A$, which
is more easily inverted since it is block upper triangular:
\begin{equation}
A
\left(
\begin{MAT}(@)[2pt]{c}
\begin{MAT}(@)[1pt,0.5cm,2.1cm]{c}
f \\
\end{MAT}\\2
\begin{MAT}(@)[1pt,0.5cm,0.4cm]{c}
\tau \\
\end{MAT}\\
\end{MAT}
\right)
\equiv
\left(
\begin{MAT}(@)[2pt]{c2c}
\begin{MAT}(@)[0pt,3.5cm,0.4cm]{c}
a\hat{e}_1^T\\ 
\end{MAT} & -Q_{N,N}\\2
\begin{MAT}(@)[0pt,3.5cm,1.5cm]{c}
L_{lo,i}\\
\end{MAT} & -Q_{lo,N} \\2
\begin{MAT}(@)[0pt,0.5cm,0.4cm]{c}
0\\
\end{MAT} & 1 \\
\end{MAT}
\right)
\left(
\begin{MAT}(@)[2pt]{c}
\begin{MAT}(@)[1pt,0.5cm,2.1cm]{c}
f \\
\end{MAT}\\2
\begin{MAT}(@)[1pt,0.5cm,0.4cm]{c}
\tau \\
\end{MAT}\\
\end{MAT}
\right)
=
\left(
\begin{MAT}(@)[2pt]{c}
(Qg)_N \\2
\begin{MAT}(@)[1pt,0.5cm,1.5cm]{c}
(Qg)_{lo} \\
\end{MAT}\\2
\beta\\
\end{MAT}
\right)
\label{eq:goodmatrix}
\end{equation}
where $a$ is an arbitrary value whose order of magnitude is that of the 
dominant values of $L$ and $\hat{e}_1^T$ is a unit row vector corresponding
to the lowest radial wavenumber present for this $m$.
The upper-left matrix in \eqref{eq:goodmatrix} is
banded and diagonally dominant and so can be stably
inverted without any pivoting.
The matrices in \eqref{eq:goodmatrix} 
and \eqref{eq:badmatrix} differ only in their first and
last rows, so their difference $v w^T$ is of rank two:
\begin{equation}
\left(
\begin{MAT}(@)[2pt]{c2r}
\begin{MAT}(@)[0pt,3.5cm,0.35cm]{c}
L_{N,i}\ -\  a\hat{e}_1^T\\ 
\end{MAT} & 0\\2
\begin{MAT}(@)[0pt,3.5cm,1.5cm]{c}
0\\
\end{MAT} & 0 \\2
\begin{MAT}(@)[0pt,0.5cm,0.35cm]{c}
B^T\\
\end{MAT} & -1 \\
\end{MAT}
\right)
=
\left(
\begin{MAT}(@)[2pt]{c}
\begin{MAT}(@)[0pt,0.75cm,0.4cm]{cc}
1 & 0 \\
\end{MAT}\\2
\begin{MAT}(@)[0pt,0.75cm,1.5cm]{cc}
0 & 0 \\
\end{MAT}\\2
\begin{MAT}(@)[0pt,0.75cm,0.4cm]{cc}
0 & 1 \\
\end{MAT}\\
\end{MAT}
\right)
\left(
\begin{MAT}(@)[2pt]{c2r}
\begin{MAT}(@)[0pt,3.5cm,0.35cm]{l}
-a\hat{e}_1^T\ +\  L_{N,i}\\ 
\end{MAT} & 0\\
\begin{MAT}(@)[0pt,0.5cm,0.35cm]{c}
B^T\\
\end{MAT} & -1 \\
\end{MAT}
\right)\equiv v w^T
\end{equation}
We define:
\begin{equation}
L'\equiv
\left(
\begin{MAT}(@)[2pt]{c}
\begin{MAT}(@)[0pt,3.5cm,0.35cm]{c}
L_{N,i}\ -\  a\hat{e}_1^T\\ 
\end{MAT}\\2
\begin{MAT}(@)[0pt,3.5cm,1.5cm]{c}
0\\
\end{MAT}\\
\end{MAT}
\right),\quad
q\equiv PQ\hat{e}_N=
\left(
\begin{MAT}(@)[2pt]{c}
\begin{MAT}(@)[0pt,0.5cm,0.35cm]{c}
Q_{N,N}\\ 
\end{MAT}\\2
\begin{MAT}(@)[0pt,0.5cm,1.5cm]{c}
Q_{lo,N}\\
\end{MAT}\\
\end{MAT}
\right)
\label{eq:Lprime}
\end{equation}
%
%\begin{subequations}
%\begin{align}
%\tau &= \beta \\
%L^\prime f & = PQg +\beta q\Longleftrightarrow 
%\left(\begin{MAT}(@)[2pt]{c}
%\begin{MAT}(@)[0pt,3.5cm,0.4cm]{l}
%a\hat{e}_1^T\\ 
%\end{MAT}\\2
%\begin{MAT}(@)[0pt,3.5cm,1.75cm]{c}
%L_{lo,i}\\
%\end{MAT}\\
%\end{MAT}\right)
%\left(
%\begin{MAT}(@)[2pt,0.5cm,2.55cm]{c}
%f\\
%\end{MAT}\right)
%=
%\left(
%\begin{MAT}(@)[2pt]{c}
%(Qg)_N \\2
%\begin{MAT}(@)[1pt,0.5cm,1.75cm]{c}
%(Qg)_{lo} \\
%\end{MAT}\\
%\end{MAT}
%\right)
%+\beta
%\left(
%\begin{MAT}(@)[2pt]{c}
%Q_{N,N} \\2
%\begin{MAT}(@)[1pt,0.5cm,1.75cm]{c}
%Q_{lo,N} \\
%\end{MAT}\\
%\end{MAT}
%\right)
%\end{align}
%\left(\begin{array}{cc} a\hat{e}_1^T & \rightarrow \\ L_{lo,i} & \rightarrow\end{array}\right) 
%\left(\begin{array}{c} f \\ \downarrow\\ \end{array}\right) =
%\left(\begin{array}{c} (Qg)_N \\ (Qg)_{lo} \end{array}\right) 
%\label{eq:Lprime}
%\end{subequations}
%
In addition to the ability to solve \eqref{eq:goodmatrix},
the Sherman-Morrison-Woodbury formula \eqref{eq:smw} requires only
the inversion of the following $2\times 2$ matrix:
\begin{align}
I+ w^T A^{-1} v &= I +
\left(
\begin{MAT}(@)[2pt]{c2r}
\begin{MAT}(@)[0pt,3.5cm,0.35cm]{l}
-a\hat{e}_1^T\ +\  \hat{e}^T_{N}L\\ 
\end{MAT} & 0\\
\begin{MAT}(@)[0pt,0.5cm,0.35cm]{c}
B^T\\
\end{MAT} & -1 \\
\end{MAT}
\right) 
A^{-1}
\left(
\begin{MAT}(@)[2pt]{c}
\begin{MAT}(@)[0pt,0.75cm,0.4cm]{cc}
1 & 0 \\
\end{MAT}\\2
\begin{MAT}(@)[0pt,0.75cm,1.5cm]{cc}
0 & 0 \\
\end{MAT}\\2
\begin{MAT}(@)[0pt,0.75cm,0.4cm]{cc}
0 & 1 \\
\end{MAT}\\
\end{MAT}
\right)
%\left(\begin{array}{lr} \hat{e}_N^T L-a \hat{e}_1^T & 0\\
%B^T & -1
%\end{array}\right) A^{-1}
%\left(\begin{array}{cc} 1 & 0 \\ 0 & 0 \\ 0 & 1\end{array}\right)
\nonumber\\ 
&=\left(\begin{array}{lr}
\hat{e}_N^T L (L^\prime)^{-1} \hat{e}_1  &  
\quad(-a \hat{e}_1^T+\hat{e}_N^T L) (L^\prime)^{-1} q \hat{e}_N\\
B^T     (L^\prime)^{-1} \hat{e}_1  &   B^T               (L^\prime)^{-1} q\hat{e}_N
\end{array}\right) &
\end{align}
where we have used the fact that $a\hat{e}^T_1 (L^\prime)^{-1} \hat{e}_1 = 1/a$.

\subsection{Nonlinear terms}
\label{sec:nonlin_term}

To compute the nonlinear terms 
\begin{subequations}
\label{eq:nonlin}
\begin{alignat}{2}
\source_\psi&\equiv&&\ez\cdot\curl\bsource\\
\source_\phi&\equiv&-&\ez\cdot\curl\curl\bsource
\end{alignat}
\end{subequations}
where
\begin{equation}
\bsource\equiv (\bu\cdot\bnabla)\bu =
\frac{1}{2}\grad(\bu\cdot\bu)-\bu\times\bom 
\label{eq:sdef}\end{equation}
we use the pseudo-spectral method
\cite{Orszag71d}, in which fields are
transformed into physical space, the nonlinear terms are carried out via
pointwise multiplication, and the results transformed back into spectral
space.  Computing the nonlinear term $\bsource$ in the rotational form 
$-\bu\times\bom$ requires only 9 spectral$\leftrightarrow$physical
transforms as compared to the 15 transforms required by the convective form
$(\bu\cdot\grad)\bu$. The difference between them
is annihilated by the curls taken in \eqref{eq:nonlin}, so we
will write $\bsource=-\bu\times\bom$.

The calculation of the nonlinear terms presents two
difficulties. The first involves radial parity and appears 
when creating $\bsource$.
We have sought to use only scalar fields which can
be represented by expansions of type \eqref{eq:fourdecomp}
and constructed using radial operators such as
\eqref{eq:parity_ops}
which preserve radial parity and can be implemented via recursion relations.
The components of velocity and vorticity, defined in cylindrical coordinates
using the toroidal and poroidal potentials as
\begin{subequations}
\begin{align}
	\bu &= \left(\ir\pdth\psi+\pdrz\phi\right)\er + \left(\ir\pdthz\phi-\pdr\psi\right)\et 
	+ \left(\lap_h\phi\right)\ez\label{eq:def_y_chap_nonlin}\\
	\bom &= \left(\pdrz\psi-\ir\pdth\lap\phi\right)\er + \left(\ir\pdthz \psi - \pdr\lap\phi\right)\et + 
	\ir\left(-\lap_h\psi\right)\ez
\end{align}
\label{eq:uom}
\end{subequations}
do not have this property. 
We therefore construct modified fields:
\begin{subequations}
\begin{align}
\bu^* &\equiv ru_r\er + ru_\theta\et + u_z\ez 
= (\pdth\psi+r\pdrz\phi)\er + (\pdthz\phi-r\pdr\psi)\et + (\lap_h\phi)\ez\\
\bom^* &\equiv r\omega_r\er + r\omega_\theta\et + \omega_z\ez 
=(\pdth u_z-\pdz u^*_\theta)\er+(\pdz u^*_r-r\pdr u_z)\et+(-\lap_h\psi)\ez
\end{align}
\label{eq:modif_u_omega}
\end{subequations}
whose components have the same parity as $\psi$ and $\phi$ and so can be
created and acted upon using the differential operators in
\eqref{eq:parity_ops}.  The modified fields $\bu^*$ and $\bom^*$ are
transformed into physical space, where their cross product is taken to form
\begin{align}
\bsource^* &\equiv r\,\source_r\er + r\,\source_\theta\et + \source_z\ez 
\equiv \source_r^*\er + \source_\theta^*\et + \source_z^*\ez =-\bu^*\times\bom^*
\end{align}

The second difficulty appears when differentiating $\bsource$ in \eqref{eq:nonlin}
and involves regularity.
A vector function which is analytic at the origin must obey
conditions analogous to \eqref{eq:polarsing}, namely
\begin{equation}
f_r(r) = r^{|m-1|}p_r(r^2) \qquad 
f_\theta(r) = r^{|m-1|}p_\theta(r^2) \qquad 
f_z(r) = r^{m}p_z(r^2)
\label{eq:vector}
\end{equation}
where $p_r$, $p_\theta$ and $p_z$ are polynomials.
We require not only regularity at the origin of $\bsource$, but also 
regularity of its curl and double curl.
When $\bsource = -\bu\times\bom$ is calculated analytically,
this is in fact the case.
However, the numerical transforms to and from physical space introduce
aliasing errors which destroy this property.  Full dealiasing would multiply
the time necessary for evaluating of the nonlinear term by a factor
$\approx4.5$.
\citet{Matsushima} suggest instead that all terms that could potentially suffer
in spectral space from singular operations (like dividing by $r$) be
evaluated in physical space (at collocation points excluding the coordinate
origin) and transformed back to the spectral space using the radial transform,
ensuring the correct polynomial order for a given Fourier mode. We have
generalized this approach to the evaluation of $\source_\psi$ and $\source_\phi$.
For each Fourier mode $m$, we write
\begin{subequations}
\begin{alignat}{6}
\source_\psi&= &&\frac{1}{r}\left(\pdr r \source_\theta - im\source_r\right)
&=&&&\frac{1}{r^2}(r\pdr -m)\source^*_\theta-\frac{im}{r^2}(\source^*_r+i\source^*_\theta)
\label{eq:nonlin_term_pot_psi}\\
\source_\phi&= &-&\frac{1}{r}\pdz\left(\pdr r \source_r + i m \source_\theta\right)+\laph \source_z
&=&-&&\frac{1}{r^2}\pdz(r\pdr -m)\source^*_r+\laph \source_z-\frac{m}{r^2}\pdz(\source^*_r+i\source^*_\theta)
\label{eq:nonlin_term_pot_phi}
\end{alignat}
\label{eq:nonlin_term_pot}
\end{subequations}
The purpose of the decomposition in \eqref{eq:nonlin_term_pot} can be
explained as follows.  Consider the first terms on the right-most-sides of
\eqref{eq:nonlin_term_pot}, those which do not contain $(\source^*_r +
i\source^*_\theta)$.  Because they do not generate terms of lower polynomial
order, these terms preserve regularity and can be carried out in spectral
space.  In contrast, although $(\source^*_r + i\source^*_\theta)$ should
be divisible by $r^2$, aliasing error in the multiplication \eqref{eq:sdef}
can generate terms of lower polynomial order.  (More specifically, $\source^*$
and $i\source^*_\theta$ are not individually divisible by $r^2$ and aliasing
perturbs the cancellation in the sum $\source^* + i\source^*_\theta$.)  This
division by $r^2$ is therefore carried out in physical space.  This increases
the number of spectral$\leftrightarrow$physical transformations only from 9 to
10. Further details can be found in \cite{BoronskiPhD}.

\subsection{Parallelization}

The separability of almost the entire algorithm (except for the nonlinear
term) between the Fourier modes makes parallelization of the code quite
straightforward. Our code was parallelized using the MPI protocol which made
it possible to run even very time-consuming three-dimensional simulations
with resolutions such as $M\times K\times N =128\times 160\times 90$.  
Spectral methods are often
considered to be poorly suited for parallelization as they require the
exchange of all the data at each timestep of the simulation.  In our code, all
the necessary data exchange is done within two calls to the
\emph{MPI\_Alltoall} MPI subroutine treating, in total, 10 three-dimensional
fields at each time step. Even though this may seem to be a large operation,
on the IBM Power4 architecture with 64 processors we found that the time
overhead per timestep due to the data exchange is counterbalanced by more
efficient usage of the processor cache memory: each processor of the
parallelized code treats smaller data portions which can more easily fit into
processor's fast internal memory (cache).  We observed that the total CPU time
used by the parallel code is often smaller than that used by the serial code
treating the same problem, as shown in table \ref{tab:parallel}.  The efficiency of the parallel code depends,
however, on the speed and the latency of the inter-processor network: the IBM
Power4 architecture has particularly fast connection between the nodes which
use the mixed model \emph{fast network/shared memory} communication between
processors. We conclude that for modern massively parallel computers,
parallelization of a pseudo-spectral code does not necessarily degrade its
efficiency but can actually enhance it.  Additional technical information
about the MPI parallelization of the code is given in \cite{BoronskiPhD}.

\begin{table}[!h]
\centering
\begin{tabular}{ccccc}
\toprule
& number of & resolution & CPU (sec)/ &
CPU (sec)/\\
& processors& $M\times K \times N$ & timestep &
(timestep$\times M$)\\
\midrule
2D & single processor & $1 \times 96 \times 48$ & 0.035 & 0.035 \\
3D & single processor & $32 \times 96\times 48$ & 1.7 sec & 0.053\\
3D & 32 processors & $32 \times 96\times 48$ & 1.5 sec & 0.047\\
\bottomrule
\end{tabular}
\caption{CPU timings on IBM Power4, as the azimuthal resolution $M$ 
and number of processors is varied.}
\label{tab:parallel}
\end{table}

\section{Tests and validation}
\label{sec:validation}

We now describe the ways in which we have validated the hydrodynamic code
described here and in our companion paper \cite{otherpaper}.  We have obtained
an exact polynomial solution to the nested Helmholtz-Poisson solver, which is
the by far most complicated portion of the code; we present its form in the
hopes it may prove useful to other researchers.  For non-polynomial solutions,
we have analyzed the effect of the corner singularity on the error and
verified the spectral convergence.  We have tested the full nonlinear
time-dependent program by simulating 2D and 3D stationary and time-periodic
rotor-stator flows which are well documented in the literature.

\subsection{Polynomial solutions}

There are no polynomial solutions which exactly satisfy the linear
differential equations \eqref{eq:potMHD_u} 
and boundary conditions \eqref{eq:bcs} with which to compare a
numerical solution.
We can, however, formulate a polynomial solution to a related problem:
the time-discretized equations 
\begin{subequations}
\begin{align}
(I-\dt\iR\lap)\lap_h\psi &= rhs_\psi\\
(I-\dt\iR\lap)\lap\lap_h\phi &= rhs_\phi
\end{align}
\end{subequations}
relating the four fields $\psi$, $\phi$, $rhs_\psi$, $rhs_\phi$
and subject to boundary conditions \eqref{eq:bcs}.
Under these conditions, no error is introduced in imposing the boundary
conditions via the $\tau$ method (see section \ref{sec:viscous}).  Comparison
between this analytic solution and a numerical solution provides a test of the
linear Helmholtz/Poisson solver, including the enforcement of the boundary
conditions via the influence matrix method, which is independent of the
evaluation of the nonlinear terms, the temporal integration and the spatial
resolution.

Using a symbolic algorithm implemented in Maple, we have
calculated polynomial solutions which contain several radial and axial basis
functions (wavenumbers) and which correspond to a realistic profile
$\angvel_\pm$, such as $\angvel_\pm^{poly}=\pm(1-r^6)$.  Polynomial solutions
were calculated for Fourier modes ranging from $m=0$ to $m=5$.  For 
$m=0$, we sought a solution containing two recirculation rolls separated
by the mid-plane $z=0$, leading to the potentials, right-hand-sides, and
velocities:
\begin{subequations}
\begin{align}
\psi^{poly}(r,z) &=\;\frac{1}{64}z(-30z^2+33z^4+5)r^8-\frac{5}{48}z(z-1)(z+1)(5z^2-1)r^6-\frac{1}{2}z^5r^2\\
rhs_\psi^{poly}(r,z)&=-z(-690z^2+33z^4+185)r^6+\frac{3}{4}z(-1970z^2+425+1609z^4)r^4\\ 
&-60z(z-1)(z+1)(5z^2-1)r^2-40z^3+2z^5\\ 
\phi^{poly}(r,z) &=-\frac{1}{2}(r-1)^3(r+1)^3(z-1)^2(z+1)^2 z\\ 
rhs_\phi^{poly}(r,z) &=-72z(5z^2-33)r^4-96z(3z^4+108-131z^2)r^2+1248z^5\\ 
& +4344z-6456z^3\\
u_r^{poly} &= -3r(z-1)(z+1)(5z^2-1)(r-1)^2(r+1)^2\\ 
u_\theta^{poly}(r,z) &=
-\frac{1}{8}zr(r-1)(r+1)(-30r^4z^2+5r^4+33r^4z^4+8r^2z^4+8z^4)\\ 
u_z^{poly}&= \; 6z(z-1)^2(z+1)^2(r-1)(r+1)(3r^2-1)
\end{align}
\label{eq:poly_sol_m0}
\end{subequations}
For all polynomial solutions tested, we obtained numerical solutions for which
the maximum relative errors in $\psi$, $\phi$ and in satisfaction of the
boundary conditions, as well as the equations, are typically
$O\left(10^{-14}\right)$, \ie~machine precision, and never exceed
$O\left(10^{-12}\right)$.

\subsection{Non-polynomial solutions: error analysis}
\label{sec:non-poly}

If the right-hand-sides $rhs_\psi$ and $rhs_\phi$ are not polynomials that can
be exactly represented within the spatial discretization, then the equations
cannot be exactly satisfied.  Note that the tau error introduced in favor of
satisfying the boundary conditions at each level of the nested system of
equations is necessarily propagated to the next level.  That is, $\psi^{n+1}$
is not an exact solution to the Poisson problem \eqref{eq:system2} and in
addition, the right-hand-side $f_\psi$ is not an exact solution to the
Helmholtz problem \eqref{eq:system1}.  This implies that the error is not
isolated in the equations corresponding to the highest wavenumbers, but is
distributed among all the equations.  However, we will see below
that non-satisfaction of the equations for low
wavenumbers does not have severe consequences.

In figures \ref{fig:eq_error_lin_m0} and \ref{fig:eq_error_lin_m2}, 
we display the relative error 
\begin{EQ}[rclcrcl]
	\epsilon_\psi &=& \frac{|E\laph\psi-rhs_\psi|}{\sup|rhs_\psi|}
	, &\qquad&
	\epsilon_\phi &=& \frac{|E\lap\laph\phi-rhs_\phi|}{\sup|rhs_\phi|}
\label{eq:rel_err_eq_def}
\end{EQ}
in physical space $(r,z)$ for the $m=0$ and $m=2$ modes.  
In this subsection,
as in the previous one, we restrict ourselves to the Stokes problem,
whose errors behave like those of the Navier-Stokes equations.  
As is usually the case, the
error is concentrated in the neighborhood of the boundaries, decaying rapidly
away from them.

\begin{psfrags}
\psfrag{r}[r][][0.75]{$r$}
\psfrag{z}[r][][0.75]{$z$}
\begin{figure}[!h]
	\centering 
\psfrag{val}[r][][0.75]{{\large$\epsilon_\psi^{m=0}$}}
	\includegraphics[width=0.45\textwidth]{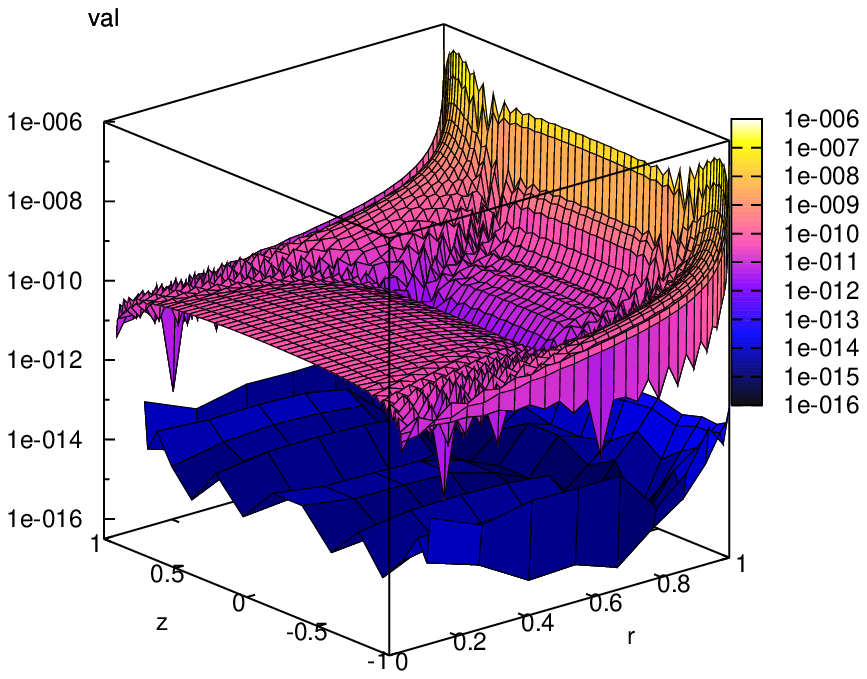}%err_psi_re500lin_m0_c.eps}
\psfrag{val}[r][][0.75]{{\large $\epsilon_\phi^{m=0}$}}
	\includegraphics[width=0.45\textwidth]{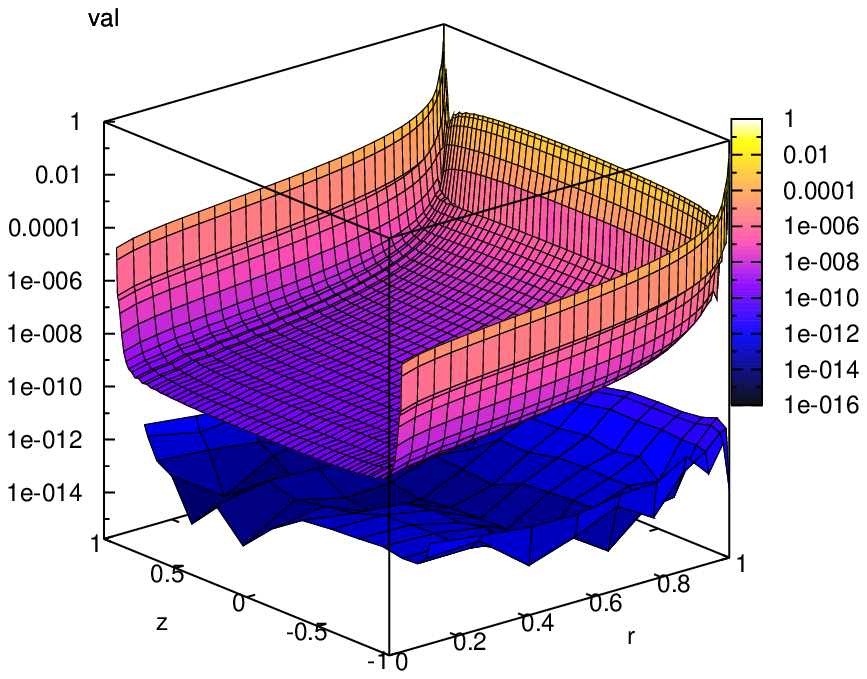}%err_phi_re500lin_m0_c.eps}
	\caption{Error of satisfaction of equations for $m=0$ after 100
	  timesteps of nested Stokes solver. Upper surfaces show the
	  error for the Stokes problem, lower surfaces for the polynomial
	  solution \eqref{eq:poly_sol_m0}.  Resolution used: $K=64,\ N=32$.
	  Left: $\epsilon_\psi^{m=0}$. Right: $\epsilon_\phi^{m=0}$.}
	\label{fig:eq_error_lin_m0}
\end{figure}
\end{psfrags}
\begin{psfrags}
\psfrag{r}[r][][0.75]{$r$}
\psfrag{z}[r][][0.75]{$z$}
\begin{figure}[!h]
	\centering
		\psfrag{val}[r][][0.75]{{\large $\epsilon^{m=2}_{eq_\psi}$}}
		\includegraphics[width=0.45\textwidth]{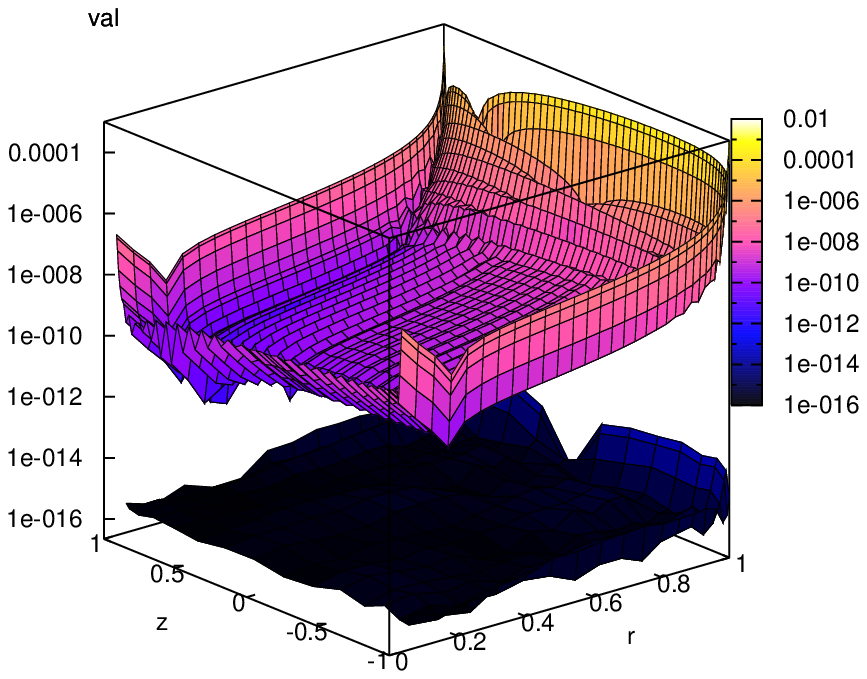}%err_psi_re500lin_m0_c.eps}
		\psfrag{val}[r][][0.75]{{\large $\epsilon^{m=2}_{eq_\phi}$}}
		\includegraphics[width=0.45\textwidth]{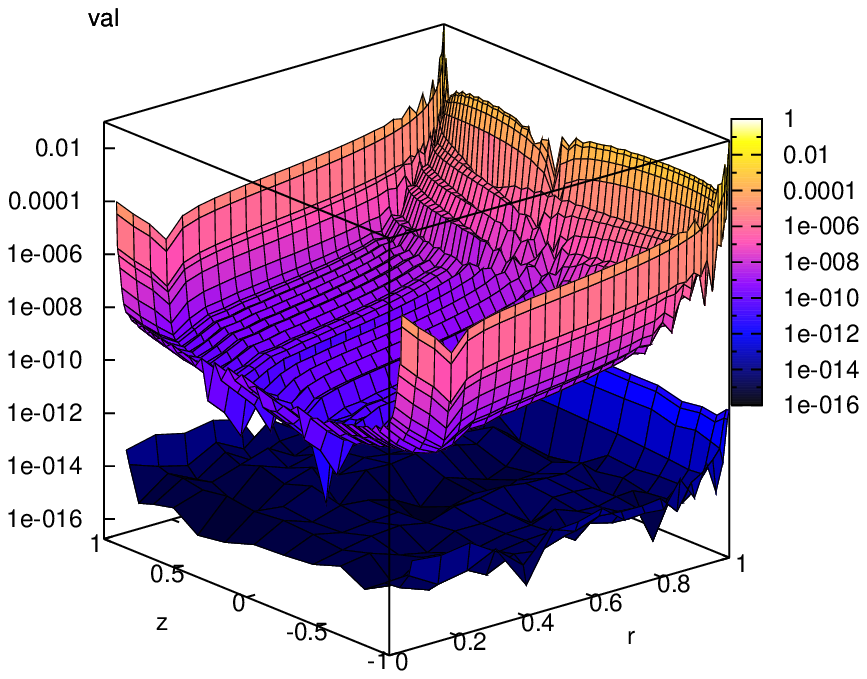}%err_phi_re500lin_m0_c.eps}
	\caption{Error of satisfaction of equations for $m=2$, after 100
	  timesteps of nested Stokes solver. Upper surfaces show the
	  error for the Stokes problem, lower surfaces for a polynomial
	  solution. Resolution used: $K=64,\ N=32$.
	  Left: $\epsilon^{m=2}_{eq_\psi}$. Right: $\epsilon^{m=2}_{eq_\phi}$.}
	\label{fig:eq_error_lin_m2}
\end{figure}
\end{psfrags}

These errors represent quite severe criteria, since
these measure the satisfaction of the curl and double curl of the original
equations; the error in satisfying the Stokes equations
themselves is considerably lower,
on the order of $\epsilon_\psi$ divided by $O(K,N)$
and $\epsilon_\phi$ divided by $O(K^2,N^2)$ near the boundary.
This estimate relies on the fact that the spectrum of the error in 
satisfying the Stokes equations is approximately uniform,
so that its derivatives are dominated by its high wavenumber terms.

The error depends significantly on the right-hand-sides $rhs_\psi$ and
$rhs_\phi$. For arbitrary right-hand-sides, the boundary conditions can be
very constraining and can lead to a nearly singular solution suffering from
spurious oscillations, with an error near the boundary that is $O(1)$.
However, $\epsilon_\psi,\epsilon_\phi$ are considerably smaller when the
right-hand-sides are calculated from the solution at the previous timestep,
especially for the Stokes problem for which the nonlinear term is zero.  

For the axisymmetric modes shown in figure \ref{fig:eq_error_lin_m0},
$\epsilon_\psi$ is $O(10^{-6})$ on the boundaries and $O(10^{-10})$
for the internal points. This is much less than $\epsilon_\phi$, which reaches
$O(1)$ at the cylinder corners.  To understand this, we recall that in the
axisymmetric case, the toroidal flow described by $\psi$ is azimuthal and the
poloidal flow described by $\phi$ is in the $(r,z)$ plane.  The azimuthal flow
described by $\psi$ follows smooth paths, while the flow described by $\phi$
must abruptly change direction near the corners.  This poloidal flow in fact
resembles the analytic asymptotic solution derived by Moffatt \cite{Moffatt}
for 2D flow in a rectangular container, which is weakly singular in
that its vorticity behaves like
$\rho^{1.74}$ (for a small distance $\rho$ from the corner).
This in turn implies that the axial component of the Laplacian 
of the Stokes equation measured by $\epsilon_\phi$
diverges for exact solutions to the continuous 2D Stokes problem,
in contrast with the numerically computed solution which 
is forced to be finite and regular.
Satisfaction of the Stokes equation itself
follows from the satisfaction of the boundary condition
which our code imposes to precision $O(10^{-14})$.

For the non-axisymmetric modes shown in figure \ref{fig:eq_error_lin_m2},
typical errors at the corners are 
$O(0.01)$ for $\psi$ and $O(0.1)$ for $\phi$ and the same corner
singularity is observed for both.
This is because the 3D Stokes solutions are also weakly singular at the 
corners \citet{Hills} and $\psi$ and $\phi$ are coupled for $m\neq 0$.
The relative error decreases rapidly away from the boundaries:
it is $O(10^{-4})$ only two gridpoints away and reaches reaches $O(10^{-6})$
for the interior points.  

The accuracy could be further enhanced -- or confined to certain modes --
by completing the method with the $\tau$-correction
\cite{Tuckerman-divfree}%,Haldenwang,Kleiser-Schumann}
, which takes into account the high-wavenumber residuals resulting
from the imposition of the boundary conditions for
each Helmholtz or Poisson problem. 

\subsection{Spectral convergence}

\begin{psfrags}
\psfrag{n}[r][][0.75]{$n$}
\psfrag{k}[r][][0.75]{$k$}
\begin{figure}[!h]
	\centering
		\psfrag{coeff}[r][][0.75]{{\large $\psi^{m=0}$}}
		\includegraphics[width=0.45\textwidth]{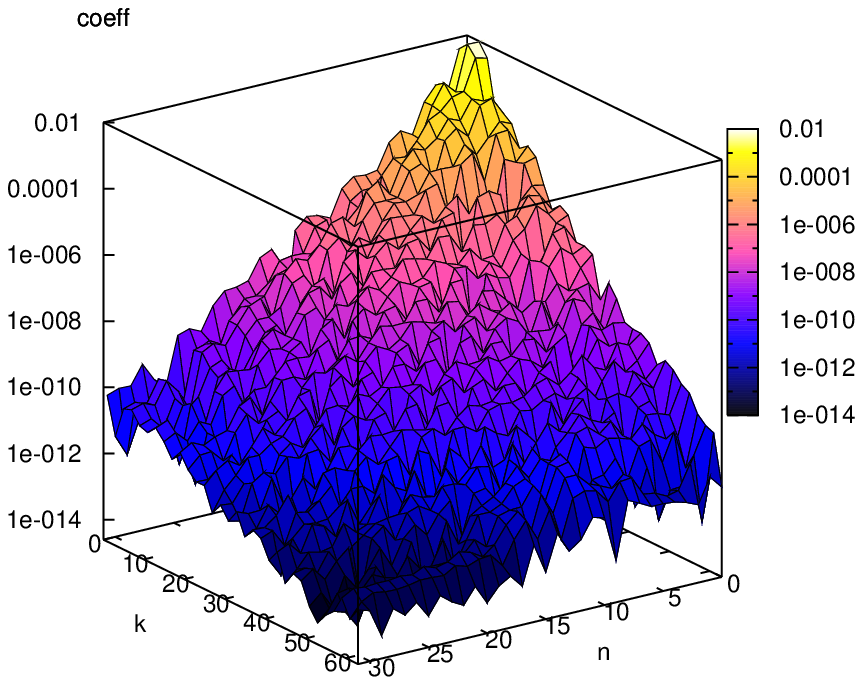}%psi_re500_m0_spec_c.eps}
		\psfrag{coeff}[r][][0.75]{{\large $\phi^{m=0}$}}
		\includegraphics[width=0.45\textwidth]{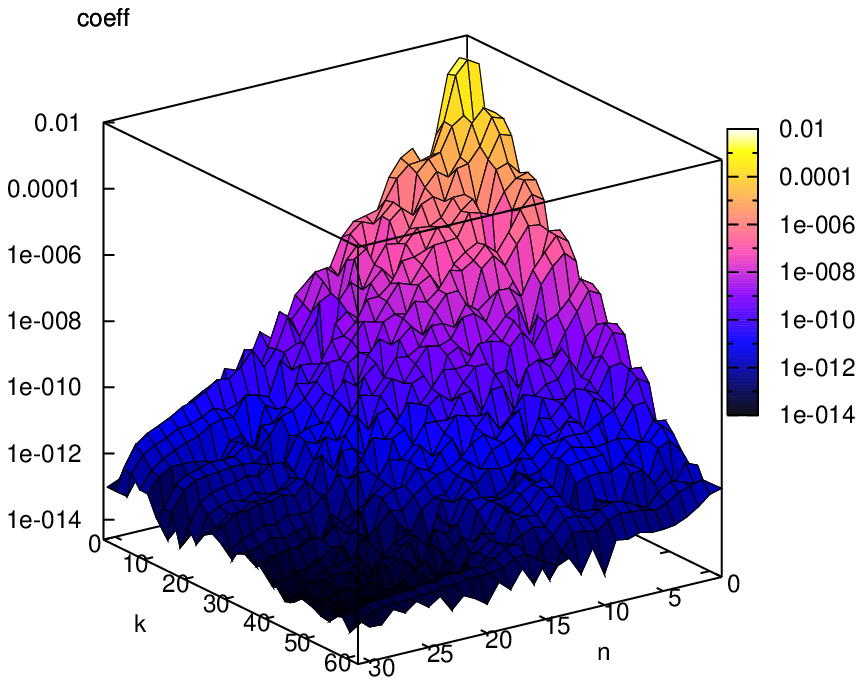}%phi_re500_m0_spec_c.eps}
	\caption{Spectral coefficients for $m=0$, after 100 timesteps of the
	nonlinear Navier-Stokes solver ($Re=750$, $\dt =0.01$). Resolution
	used: $K=64,\ N=32$. Left: $psi^{m=0}$ Right: $phi^{m=0}$.}
	\label{fig:converge_spectrum_m0_nonlin}
\end{figure}
\end{psfrags}

The important indicator of spatial convergence for spectral methods is the
decay rate of high-wavenumber coefficients in the solution fields. For a
well-behaved solver, in the absence of volume and boundary singularities, the
magnitude of spectral coefficients should decay rapidly with wavenumber.  For
a laminar flow this decay rate should be exponential, but the presence of thin
boundary layers can significantly influence the convergence.  Figure
\ref{fig:converge_spectrum_m0_nonlin} shows the spectral convergence for a
full Navier-Stokes simulation at $Re=750$ at $T=100\dt=1$ from initial
conditions of $u_\theta=\omega_+(r) 2z/h, u_z=u_r=0$.  We show the $(r,z)$
spectral coefficients for Fourier mode $m=0$; those for $m\neq 0$ are similar.
The convergence of the spectra can be qualified as quasi-exponential, meaning
that the high-wavenumber spectral coefficients seem not to decrease below a
level $O(10^{-12})$.  This can almost certainly be attributed to the singular
character of the solution to the Stokes equation near the cylinder corners.
In support of this explanation, we note that the spectrum obtained after 100
timesteps of the linear (Stokes) solver behaves similarly, except for
$\psi^{m=0}$, which displays true explonential decay, down to levels of
$10^{-22}$ for the same resolution.

\subsection{Axisymmetric rotor-stator configuration}
\begin{figure}[h]
\centering
\includegraphics[width=0.372\textwidth]{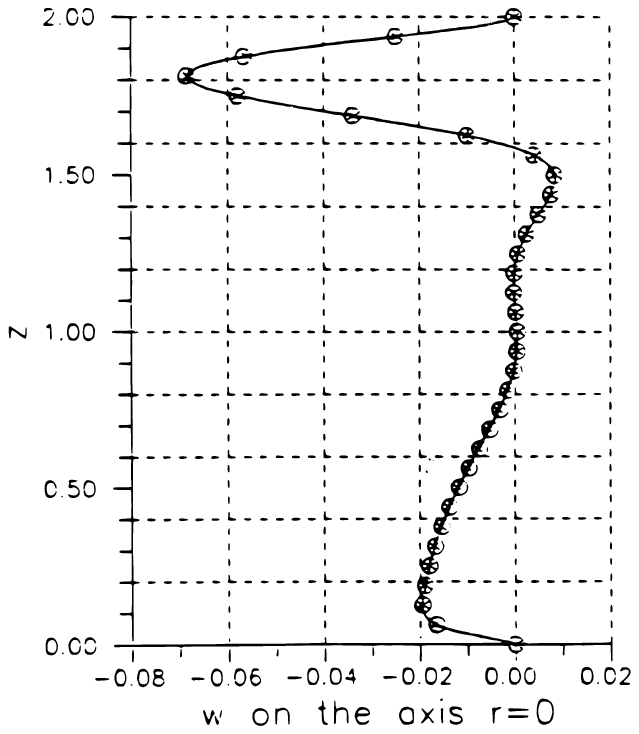}\quad
\includegraphics[width=0.36\textwidth]{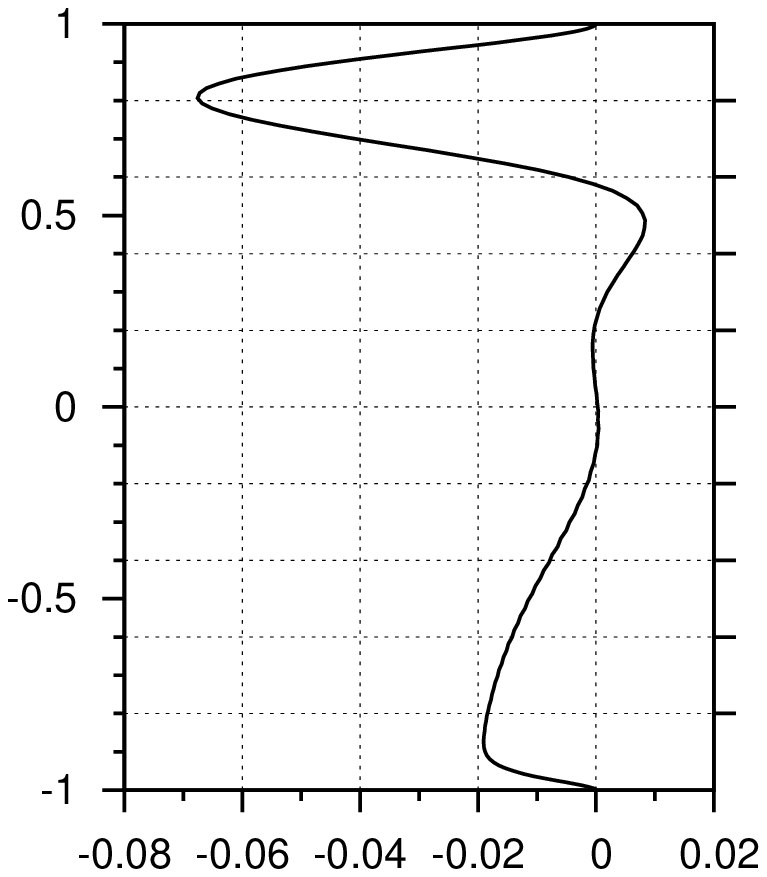}
\caption{Profile of $u_z$ at $r=0$ for rotor-stator configuration at
$Re=1850$, $h=2$.  a) Profile from \citet{Daube92}: results obtained using
$\eta-\omega$ (*) and $\bu-\omega$ ($\circ$) are superposed. b) Profile at
$t=3000$ obtained from the present code, using poloidal-toroidal decomposition
$\psi-\phi$.}
\label{fig:axial_vel}
\end{figure}

The code was first tested on the well-documented axisymmetric 
rotor-stator configuration with aspect ratio $h=2$. The first test is the
reproduction of the characteristic steady state for $Re\approx1850$ where the
flow exhibits two recirculation bubbles (one large and the other much smaller)
situated approximately at $(r=0,\,z=1/2)$ and $(r=0,\,z=0)$.
\begin{figure}[!h]
\begin{minipage}[l]{5cm}
%\vspace{-0.6cm}
\reflectbox{\includegraphics[height=8.6cm]{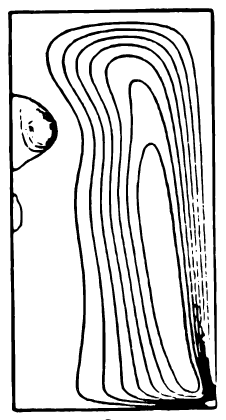}}
\caption{Contours of Stokes poloidal streamfunction $\eta$ from \cite{Daube92}
for $Re=1850$, $h=2$.}
\label{fig:daube_Psi}
\end{minipage}
\quad
\begin{minipage}[l]{10cm}
\vspace{0.85mm}
\scalebox{1.04}[1.04]{\includegraphics[width=8cm]{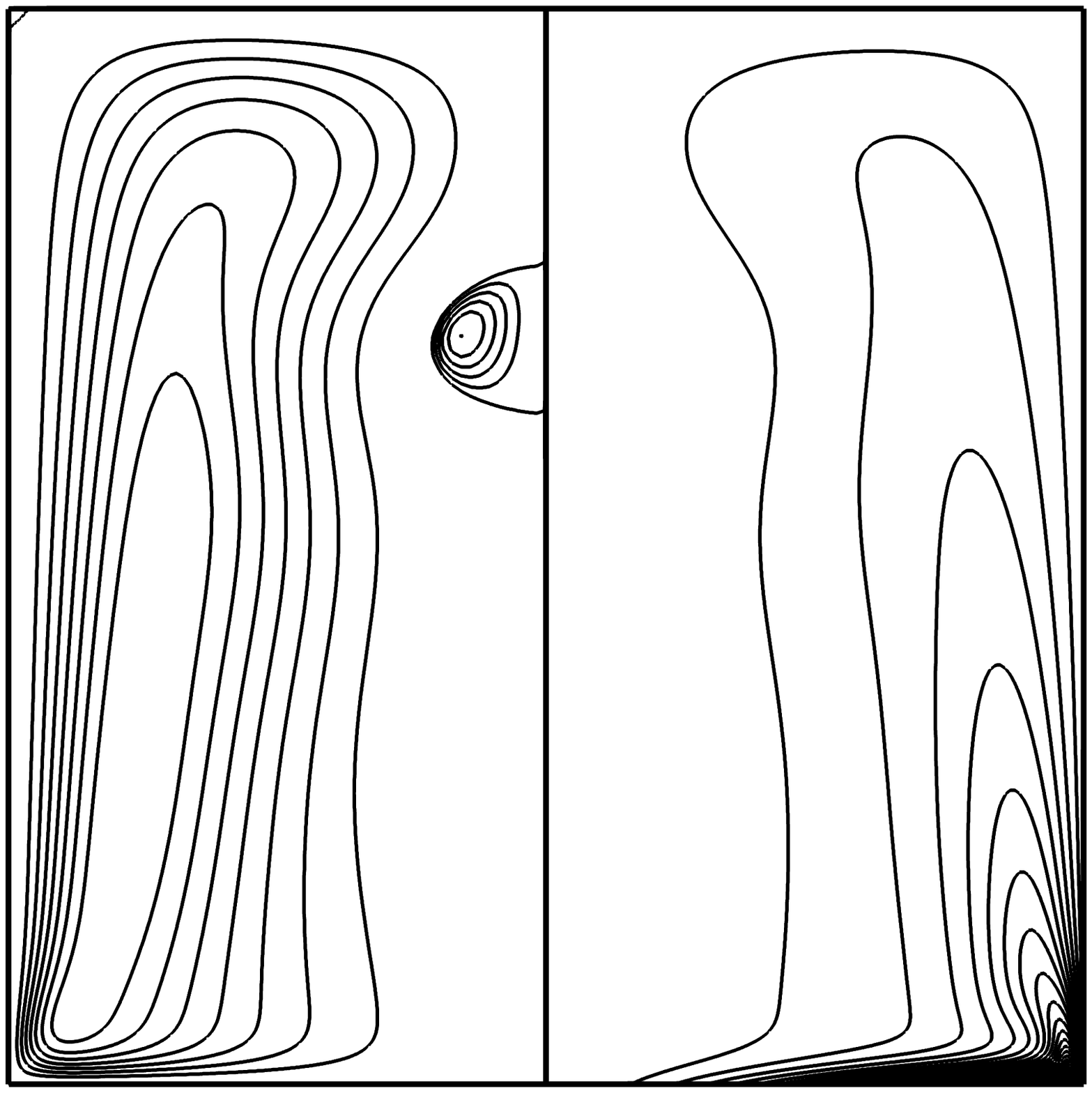}}
\vspace{1.9mm}
\caption{Rotor-stator configuration for $Re=1850$ and $h=2$.  Left: poloidal
streamfunction $\eta$. Right: toroidal streamfunction $\sigma$.}
\label{fig:rotor-stator}
\end{minipage}
\end{figure}
The contour plot of the Stokes streamfunctions defined as 
\begin{equation}
\sigma(r,z)=-r\pdr\psi(r,z)\qquad\eta(r,z)=-r\pdr\phi(r,z)
\end{equation}
presented on figure \ref{fig:rotor-stator} matches that presented
by \citet{Daube92} (figure \ref{fig:daube_Psi}) and is similar to those of
\citet{LopezShen1998} obtained for the slightly larger aspect ratio $h=2.5$.  
Quantitative agreement between our results and the previous calculations by
\citet{Daube92} and \citet{Lugt} is established by comparing the profiles of
axial velocity $u_z$ on the cylinder axis (see figure
\ref{fig:axial_vel}). This test shows excellent agreement between our
results obtained using the poloidal-toroidal formulation (fig
\ref{fig:axial_vel}b) and the velocity-vorticity formulation
(fig. \ref{fig:axial_vel}a).

\begin{figure}[h]
\begin{minipage}[lt]{0.48\textwidth}
\includegraphics[width=\textwidth,trim=0 0.3cm 0 0.2cm]{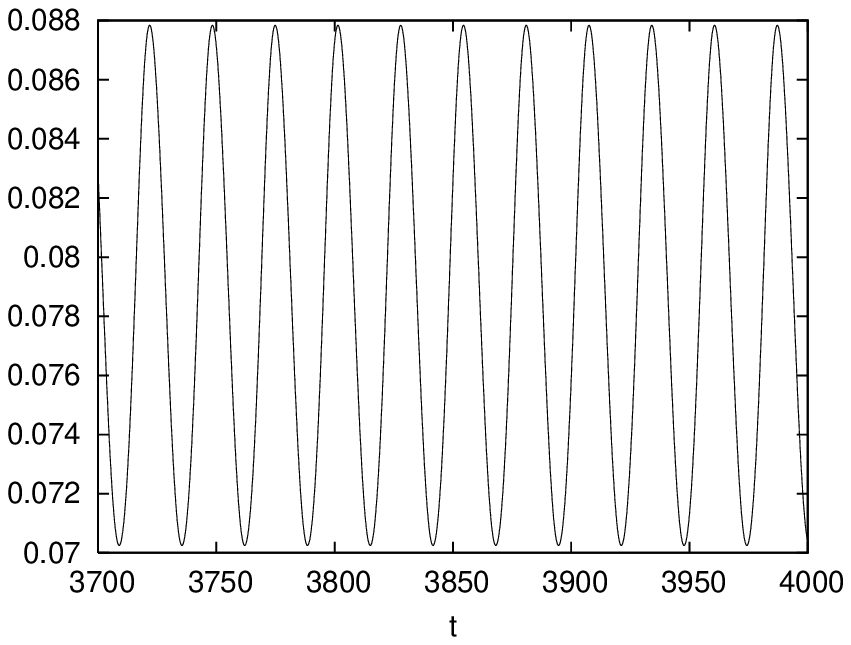}
\caption{Saturated state of time evolution of $u_\theta(r=0.5,z=0)$ for the
rotor-stator configuration at $Re=2800$, $h=2$.}
\label{fig:oscil}
\end{minipage}
\quad
\begin{minipage}[lt]{0.48\textwidth}
\includegraphics[width=\textwidth]{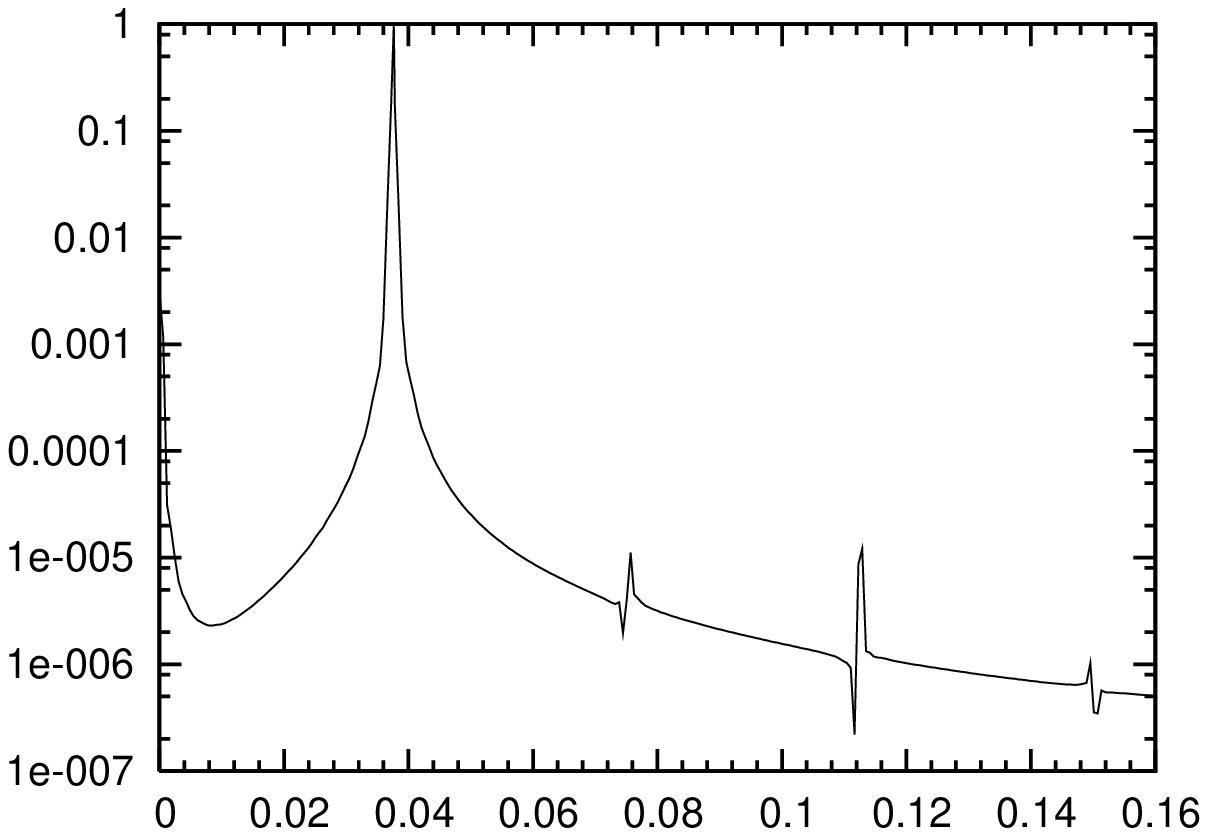}
\caption{Density power spectrum of $u_\theta(r=0.5,z=0)$ from figure \ref{fig:oscil}.}
\label{fig:oscil_spec}
\end{minipage}
\end{figure}

It was observed experimentally by \citet{Escudier} and numerically by
\citet{DaubeSorensen89}, \citet{Lopez90}, \citet{Daube92}, \citet{Gelfgat} and
\citet{Speetjens} that this flow undergoes a Hopf bifurcation toward a flow
oscillating with a frequency approximately 0.25 times the rotation
frequency. This transition occurs at $Re$ near 2600 with a period of
$T=26.55$. These values are in the ranges previously found for
this configuration (see table \ref{tab:periods}); deviations
can probably be attributed to the differences in the treatment of the
boundary conditions (regularization). The simulation was performed with
$\dt=0.01$ and with high spatial resolution $K\times N = 140\times70$ in
order to well represent the sharp regularization profile corresponding to
$\delta=0.06$ imposed on the lateral boundary $r=1$ (see section
\ref{sec:cornersing}) as proposed by \citet{LopezShen1998} and also used
by \citet{Speetjens}.
The time evolution $u_\theta(r=0.5,z=0,t)$, along with the normalized
power spectrum, are shown on figures \ref{fig:oscil} and \ref{fig:oscil_spec}.

\parbox{0.48\linewidth}{
\begin{tabular}[h]{@{}llc@{}}
\toprule
Method & Reference & $T$ \\
\midrule
$\bu-\omega$ & \citet{Daube92}  & $25.52$ \\
$\eta-\omega$ & \citet{Daube92}  & $25.84$ \\
$\bu-\omega$ & \citet{Speetjens}	& $26.61$ \\
$\bu-p$ & \citet{Gelfgat}	& $\approx 26.7$ \\
$\psi-\phi$ & this work	& $26.55$ \\
\bottomrule
\end{tabular}
%\vspace{3mm}
\captionof{table}{Oscillation period in rotor-stator configuration with $h=2$ and
$Re=2800$.}
\label{tab:periods}
}
\quad
\parbox{0.48\linewidth}{
\begin{center}
\includegraphics[width=3.75cm]{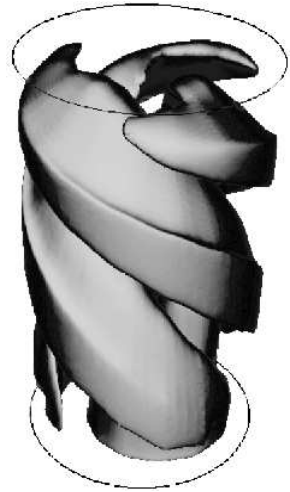}
\vspace{-1mm}
\captionof{figure}{Rotor-stator configuration with $Re_c\approx2150$ and $h=3.5$.
Isosurface $u_z\approx 0$ of axial velocity field after subtraction
of axisymmetric component.
}
\end{center}
\label{fig:m3critical}
}

\begin{figure}[h]
\centering
\begin{minipage}[lt]{0.49\textwidth}
\includegraphics[width=7cm,trim=0 -0.3cm 0 0 ]{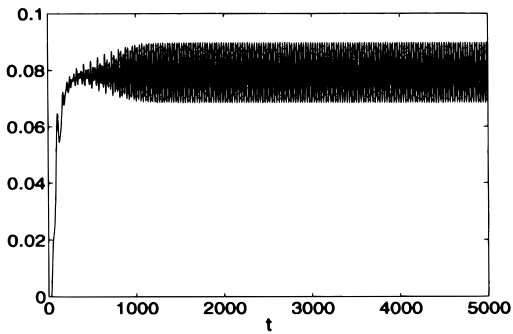}
\end{minipage}
\begin{minipage}[lt]{0.49\textwidth}
\includegraphics[width=7.2cm,height=4.8cm]{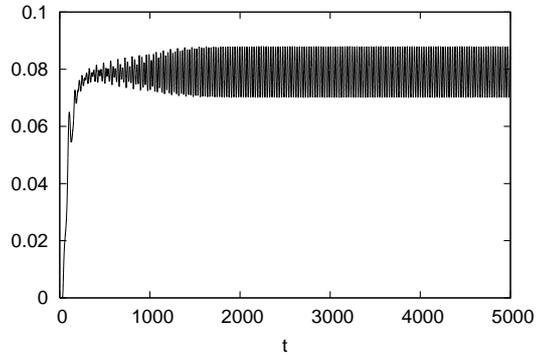}
\end{minipage}
\caption{Time history of $u_\theta(r=0.5,z=0)$ for the rotor-stator configuration at $Re=2800$, $h=2$. Left: from \citet{Speetjens}. Right: current work.}
\label{fig:time_evol}
\end{figure}
\subsection{First 3D instability\label{sec:test3d}}
We have tested the non-axisymmetric aspects of the code
using another rotor-stator configuration. For the aspect ratio of $h=3.5$ we
found that the first bifurcating mode has wavenumber $m=3$ and critical
Reynolds number $Re_c=2116$. This result is in good agreement with
\citet{Gelfgat}, where the critical Reynolds number was estimated at
$Re=2131$. A characteristic spiral analogous to that visualized by
\citet{LopezMarquesShen2001} for the same configuration is represented on
figure \ref{fig:m3critical}.

\section{Conclusion}

Motivated by a need for a numerical tool adequate for investigating
cylindrical von K\'arm\'an flow, we have written a spectral code which
solves the Navier-Stokes equation for
an incompressible fluid in a finite cylinder.
This task encompasses a number of algorithmic challenges,
which we list in order of decreasing difficulty.
The first is to impose incompressibility, a goal 
which we achieved by using the poloidal-toroidal decomposition.
In a finite cylinder, this formulation leads to 
differential equations of higher order and coupled boundary conditions
\cite{Marques90,Marques93}. In a companion paper \cite{otherpaper}, 
we described the way in which we used the influence matrix technique to
transform these equations and boundary conditions into an equivalent set of
decoupled Helmholtz or Poisson problems, each with Dirichlet boundary
conditions.

The second challenge, which has been the main focus of this article, 
is the treatment of the cylindrical axis.
The singularity engendered by the use of cylindrical coordinates
is only apparent and should not be transmitted to the fields.
We have dealt with the axis singularity by using a basis of 
radial polynomials developed by \citet{Matsushima} which are analytic at the axis
and which have properties similar to those of the Legendre polynomials. 
We have extended \cite{Matsushima} in several ways.
(i) The radial basis was developed for two dimensions, 
\ie polar coordinates. Here, we have used it to represent fields 
in a three-dimensional cylinder of finite (non-periodic) axial length.
This fairly straightforward extension was implemented 
by diagonalizing the differential operators in the axial direction.
(ii) Numerical inaccuracy in the pseudo-spectral evaluation of the 
nonlinear term can generate terms which are not analytic.
We have generalized to the high-order equations resulting from 
the poloidal-toroidal decomposition the procedure
developed by \citet{Matsushima} to avoid this problem.
Care must be taken to preserve order and parity at each stage 
of the calculation.
(iii) Differential operators expressed in this basis can,
as in most such cases \cite{Tuckerman-banded},
be reformulated as recursion relations, which can be
used to reduce the time for action or inversion of
differential operators.
Using the Sherman-Morrison-Woodbury formula, we have
formalized the stable algorithm given in \cite{Matsushima}
for solving Poisson and Helmholtz problems in a time 
proportional to the number or radial gridpoints or modes.

The third challenge is the genuine singularity at the corners
of the domain, where the disks and the cylinder which bound
the domain meet. Finite difference and finite element
codes intrinsically smooth the singularity; in contrast,
spectral expansions attempt to converge to the discontinuity,
leading to spurious oscillations.
In our code, we replaced the discontinuous boundary conditions
by a profile on the disks which is steep but continuous as
the corner is approached.
The geometrical singularity remains, but it is weak and
does not prevent spectral convergence.

Tests performed for analytic polynomial solutions to the Helmholtz
problem with an appropriate right-hand side showed that the solver
reproduces exact solutions to nearly machine precision. For
a non-polynomial solution, the solver displays
exponential convergence of spectral coefficients of the solution. 
The potential equations (corresponding to the curl and double
curl of the Navier-Stokes equations) are satisfied, with an
error of only $O(10^{-10})$ for the interior points. 
The numerical code was parallelized using the MPI protocol. This made it
possible to simulate nearly turbulent flow for $Re=5000$ with a spatial
resolution of $128\times160\times90$.

Finally, we have validated the hydrodynamic code by testing it against
well-documented problems in the literature, demonstrating
the feasibility of calculating solutions to the time-dependent
Navier-Stokes equations in a finite cylindrical geometry which
are both analytic and divergence-free.

\section*{Acknowledgements}
The simulations reported in this work were performed on the
IBM Power 4 of the Institut de D\'eveloppement et des Ressources en
Informatique Scientifique (IDRIS) of the Centre National de la Recherche
Scientifique (CNRS) as part of project 1119.
P.B. was supported by a doctoral and post-doctoral grant from the
Ecole Doctorale of the Ecole Polytechnique and the Minist\`ere de l'Education
Nationale, de la Recherche et de la Technologie (MENRT).

\bibliography{bibliography}
\bibliographystyle{elsart-num-names}
\end{document}